\documentclass[sn-apa]{sn-jnl}

\usepackage{amsthm}
\usepackage{amsmath}
\usepackage{amssymb}
\usepackage{bbm}
\usepackage{enumerate}

\jyear{2021}%

\theoremstyle{thmstyleone}%
\newtheorem{theorem}{Theorem}
\newtheorem{lemma}{Lemma}%

\theoremstyle{thmstyletwo}%
\newtheorem{remark}{Remark}%
\newtheorem{corollary}{Corollary}%

\theoremstyle{thmstylethree}%

\begin{document}

\title[Global-Local Shrinkage Priors for Asymptotic Point and Interval...]{Global-Local Shrinkage Priors for Asymptotic Point and Interval Estimation of Normal Means under Sparsity}


\author*[1]{\fnm{Zikun} \sur{Qin}}\email{qinzikun@ufl.edu}

\author[1]{\fnm{Malay} \sur{Ghosh}}\email{ghoshm@ufl.edu}

\affil[1]{\orgdiv{Department of Statistics}, \orgname{University of Florida}, \orgaddress{\street{102 Griffin-Floyd Hall}, \city{Gainesville}, \postcode{32611}, \state{Florida}, \country{U.S.A.}}}


\abstract{The paper addresses asymptotic estimation of normal means under sparsity.
The primary focus is estimation of multivariate normal means where we obtain exact asymptotic minimax error under global-local shrinkage prior.
This extends the corresponding univariate work of \citet{Prasenjit2017}.
In addition, we obtain similar results for the Dirichlet-Laplace prior as considered in \citet{Bhattacharya2015}.
Also, following \citet{vanDerPas2017}, we have been able to derive credible sets for multivariate normal means under global-local priors.}




\maketitle

\section{Introduction}\label{sec1}

Estimation of normal means under sparsity started a while ago.
Its importance is felt for the analysis of high dimensional data. 
For example, in microarray experiments, there is a multitude of genes, but only a few have impact on
a certain disease.
A foundational article appears in \citet{Donoho1992}, who provided an asymptotic 
minimax estimation rate for estimation of normal means with a large majority of zeros,
but with also a few significant departures from zeros.
The idea was pursed in a Bayesian framework by \citet{Castillo2012}
who provided the same asymptotic minimax estimation rate under a class of priors $\Pi_n$
with ``exponential decay'' (see (2.2) of their paper for its definition).

The present work addresses the same problem, but stems from several recent excellent articles
of \citet{Bhattacharya2015}, \citet{vanderPas2014},
\citet{vanderPas2016} and \citet{Prasenjit2017}.
In particular, our paper has more direct structural connection with the work of 
\citet{Prasenjit2017}, but extends their work in certain directions.

It may be pointed out that the priors of \citet{Bhattacharya2015} or \citet{Prasenjit2017} can 
be brought under the general framework of \citet{vanderPas2014}, 
but each has its own salient features which enable one to provide a more concrete set of results.
In particular, these priors, now commonly referred to as ``global-local'' priors following 
\citet{Carvalho2009, Carvalho2010}, are scale mixtures of normal priors
with the scale parameters involving both global and local components. The global components try to shrink the normal means towards 
zero,
while the local parameters try to balance the same with the end of identifying and distinguishing the true signals from the noises.
While the work of \citet{Prasenjit2017} considers a single global parameter and utilizes the same as a tuning parameter, \citet{Bhattacharya2015} considered essentially multiple global parameters and assigned certain priors on them.
These ideas will be made more specific in the following sections.

We first find the asymptotic minimax error for estimation of multivariate normal means under sparsity in the nearly-black sense \citep{Castillo2012}. It is the same as the asymptotic minimax error in the univariate case, which was proved by \citet{Donoho1992}.

We then consider estimation of multivariate normal means under global-local priors.
Like \citet{Prasenjit2017}, we obtain exact asymptotic minimaxtity results as well in this situation.
Further, in the framework of \citet{Bhattacharya2015},
we obtain asymptotic minimaxity results in the multivariate case. This is the case where we put priors to the global parameters.

The final contribution of this paper is finding credible sets for multivariate normal means following the framework of \citet{vanDerPas2017}, who considered the univariate case.
We have considered a general class of global-local priors which includes the now famous horseshoe prior,  as well as a more specific class of priors which is in the framework of \citet{Bhattacharya2015}.
Like \citet{vanDerPas2017}, we have been able to identify parameter vectors for which the posteriors give good coverage, and others for which they do not.

The outline of the remaining sections is as follows.
In Section \ref{sec:minimaxerror}, we find the asymptotic minimax error in the multivariate setting.
In Section \ref{sec:estimation}, we consider estimation of multivariate normal means and obtain exact asymptotic minimax error. We also find out the corresponding posterior contraction rates around both the estimator and the true means. 
Section \ref{sec:sets} addresses results related to credible sets of multivariate normal means.
Some final remarks are made in Section \ref{sec:remark}. The proofs of some technical lemmas are given in Appendix \ref{sec:appendixa}. The proofs of the main theorems are given in Appendix \ref{sec:appendixb}.

\section{Point Estimation of Multivariate Normal Means}
\label{sec:multivariate}

\subsection{Asymptotic Minimax Error under Nearly-Black Sparsity}\label{sec:minimaxerror}
Suppose $X_i \overset{ind}{\sim} N(\theta_i, 1), i = 1, \dots, n$.
To estimate multiple normal means, \citet{Prasenjit2017}
used a general global-local prior in which the global parameter is treated as a tuning parameter.
The exact asymptotic minimaxity was established under the prior. 
There, the true means ${\theta}_{0i}$ $(1 \leq i \leq n)$ are assumed to be sparse in the nearly-black sense \citep{Donoho1992, Castillo2012}, 
meaning that the cardinality of the non-zero ${\theta}_{0i}$'s, say $q_n$,
is $o(n)$, as $n \rightarrow \infty$.
The set of nearly-black mean vectors is denoted by
$l_0[q_n] = \{\boldsymbol{\theta} \in \mathbb{R}^n: \sum_{i=1}^n \mathbbm{1}(\theta_{i} \neq 0) \leq q_n\}$ with $q_n = o(n)$.
\citet{Donoho1992} provides the asymptotic minimax error,
\begin{equation}
    \inf_{\widehat{\boldsymbol{\theta}}} \sup_{\boldsymbol{\theta}_0 \in l_0[q_n]}
    \sum_{i=1}^n E_{\theta_{0i}} \left(\widehat{\theta} - \theta_{0i}\right)^2 = 2 q_n log\left(\frac{n}{q_n}\right)(1+o(1)), \text{ as } n \rightarrow \infty. 
\end{equation}

In the multivariate situation, 
the true means $\boldsymbol{\theta}_{0i}$ $(1 \leq i \leq n)$ being assumed to be sparse in the nearly-black sense also means that $\sum_{i=1}^n \mathbbm{1}(\boldsymbol{\theta}_{0i} \neq \boldsymbol{0}) \leq q_n$ with $q_n = o(n)$.
We denote the set of nearly-black multivariate means by
$L_0[q_n] = \{ \{\boldsymbol{\theta}_{0i}\}_{i=1}^{n} : \boldsymbol{\theta}_{0i} \in \mathbb{R}^k, i=1,\dots,n, \sum_{i=1}^n \mathbbm{1}(\boldsymbol{\theta}_{0i} \neq \boldsymbol{0}) \leq q_n\}$.
One can prove that, in the multivariate setting, the asymptotic minimax error using the Mahalanobis distance loss is the same as the asymptotic minimax error using the squared error loss in the univariate setting.
We use $\lVert \cdot \rVert_{\Sigma}$ to denote the Mahalanobis norm, e.g., $\lVert \boldsymbol{X}_i \rVert_{\Sigma}^2 = \boldsymbol{X}_i^T\boldsymbol{\Sigma}^{-1}\boldsymbol{X}_i$, where $\boldsymbol{\Sigma}$ is the positive definite population covariance matrix. 

\begin{theorem}\label{thm:minimax.error}
Suppose that $\boldsymbol{X}_i \sim N_k(\boldsymbol{\theta}_i, \boldsymbol{\Sigma})$, independently, for $i=1,\dots,n$, with a positive definite covariance matrix $\boldsymbol{\Sigma}$, and that the true mean vectors $\{\boldsymbol{\theta}_i\}_{i=1}^n$ are sparse in the nearly-black sense. If we measure the error of an estimator using the Mahalanobis distance loss, then, as $n \rightarrow \infty$, 
    \begin{equation}\label{minimaxrate}
    \inf_{\{\widehat{\boldsymbol{\theta}_i}\}} \sup_{\{\boldsymbol{\theta}_{0i}\} \in L_0[q_n]}
    \sum_{i=1}^n E_{\boldsymbol{\theta}_{0i}} \lVert\widehat{\boldsymbol{\theta}}_i - \boldsymbol{\theta}_{0i}\rVert_{\Sigma}^2  = 2 q_n log\left(\frac{n}{q_n}\right)(1+o(1)). 
\end{equation}
\end{theorem}

\begin{remark}\label{rem:norm}
    When $\boldsymbol{\Sigma}$ is positive definite, the Mahalanobis norm $\lVert \cdot \rVert_{\Sigma}$ is equivalent to the $l_2$-norm $\lVert \cdot \rVert_2$ in the sense of equivalent norms, which means that there exist two positive real constants $c$ and $C$ such that $0 < c \leq C$, for each $\boldsymbol{X} \in \mathbb{R}^k$, $c \lVert \boldsymbol{X} \rVert_2 \leq \lVert \boldsymbol{X} \rVert_{\Sigma} \leq C \lVert \boldsymbol{X} \rVert_2$. Specifically, $c = \lambda_{\max}^{-1}(\boldsymbol{\Sigma})$, inverse of the largest eigenvalue of $\boldsymbol{\Sigma}$, and $C = \lambda_{\min}^{-1}(\boldsymbol{\Sigma})$, inverse of the smallest eigenvalue of $\boldsymbol{\Sigma}$.
    So, Theorem \ref{thm:minimax.error} will not give us an exact asymptotic minimax error under the $l_2$-norm unless $\boldsymbol{\Sigma}$ satisfies certain eigenvalue conditions. Instead, we can get both lower and upper bounds of the minimax error under the $l_2$-norm. Since both bounds are of the same rate, $2 q_n log\left(n/q_n\right)(1+o(1))$, the minimax error under the $l_2$-norm must be of the same rate as well, and will only differ from it by up to a constant factor.
\end{remark}

\subsection{Minimax Estimation of Multivariate Normal Means}\label{sec:estimation}
Now, we first extend the results of \citet{Prasenjit2017} to the multivariate case.
We begin with a general global-local prior model
\begin{itemize}
    \item[(i)] $\boldsymbol{X}_{i} \vert \boldsymbol{\theta}_{i} \overset{ind}{\sim} N_k(\boldsymbol{\theta}_{i}, \boldsymbol{\Sigma})$, $i = 1, \dots,n $, 
    $\boldsymbol{\Sigma}$ is known positive definite;
    \item[(ii)] $\boldsymbol{\theta}_{i} \vert \lambda_i^2 \overset{ind}{\sim} N_k(\boldsymbol{0}, \lambda_i^2 \tau_n\boldsymbol{\Sigma}) $, $i = 1, \dots,n $, 
    where $\tau_n \in (0,1)$ is a sequence of positive constants to be chosen later, $\tau_n \rightarrow 0$ as $n \rightarrow \infty$;
    \item[(iii)] $\pi (\lambda_i^2) = K (\lambda_i^2) ^{-a-1} L(\lambda_i^2)$, $i = 1, \dots,n $, where $a > 0$ and $L$ is a slowly varying function.
\end{itemize}
In this model, the global parameter $\tau_n$ is assumed to be a tuning parameter.
Note that the horseshoe prior \citep{Carvalho2009, Carvalho2010} is a special case of this prior in
the univariate setup with $a=1/2$.

The following regularity assumptions are made:
\begin{itemize}
    \item[(I)] L is non-decreasing in its argument with $0 < m \leq L(u) \leq M < \infty$;
    \item[(II)] $0 < \lambda_{\min} (\boldsymbol{\Sigma}) \leq \lambda_{\max} (\boldsymbol{\Sigma}) < \infty$,
    where $\lambda_{\min} (\boldsymbol{\Sigma})$ and $\lambda_{\max} (\boldsymbol{\Sigma})$ denote the minimum and maximum eigenvalues of $\boldsymbol{\Sigma}$.
\end{itemize}

We estimate $\boldsymbol{\theta}_i$ using the posterior means under the global-local prior, i.e.,
\begin{equation}
    \widehat{\boldsymbol{\theta}}_i = E(\boldsymbol{\theta}_i  \mid  \boldsymbol{X}_i) = E(1-\kappa_i  \mid  \boldsymbol{X}_i)\boldsymbol{X}_i,
\text{ where } \kappa_i = (1+\lambda_i^2 \tau_n)^{-1},
\end{equation}
and $\kappa_i$ is the shrinkage factor.
The estimators using a prior (iii) are denoted by $\widehat{\boldsymbol{\theta}}^R_i$ specifically.

We prove the following theorem under this model, in which \eqref{thm1} and \eqref{thm2} concern the error contributed by the zero and non-zero true means, respectively, and \eqref{thm3} is then immediate following the previous two results. In particular, when $0 < a \le 1$, we already have an upper bound for the error. Theorem \ref{thm:minimax.error} provides the minimax lower bound, which matches the upper bound here. This fact actually finishes both the proofs of Theorem \ref{thm:minimax.error} and \eqref{thm3.second}. As shown in Theorem \ref{thm:GL}, this general class of global-local priors attains the asymptotic minimax rate in the multivariate setting, and when $0 < a \le 1$, it attains the exact asymptotic minimax error.

\begin{theorem}\label{thm:GL}
Assume that the true means are sparse in the nearly-black sense.
Under the regularity assumptions (I) and (II), using the global-local prior with a tuning parameter, i.e., a model satisfying (i), (ii) and (iii), if $\tau_n = \left(q_n/n\right)^{\frac{1+\epsilon_n}{a}}$, where $\epsilon_n = 1 / \log(\log(n/q_n))$, then,
\begin{equation}\label{thm1}
\lim\limits_{n \rightarrow \infty} { \sup_{\{\boldsymbol{\theta}_{0i}\} \in L_0[q_n]}} {\sum\limits_{i: \boldsymbol{\theta}_{0i} = \boldsymbol{0}} E_0 \lVert \widehat{\boldsymbol{\theta}}^R_i \rVert_{\Sigma}^2} /
{\left( q_n \log\left(\frac{n}{q_n}\right) \right)} = 0.
\end{equation}
and
\begin{equation}\label{thm2}
\limsup\limits_{n \rightarrow \infty} { \sup_{\{\boldsymbol{\theta}_{0i}\} \in L_0[q_n]}}
\frac{\sum\limits_{i: \boldsymbol{\theta}_{0i} \neq \boldsymbol{0}} E_{\boldsymbol{\theta}_{0i}} \lVert\widehat{\boldsymbol{\theta}}^R_i - \boldsymbol{\theta}_{0i}\rVert_{\Sigma}^2} 
{2 q_n \log(n/q_n)} \leq 1.
\end{equation}
Consequently,
\begin{equation}\label{thm3}
\limsup\limits_{n \rightarrow \infty} { \sup_{\{\boldsymbol{\theta}_{0i}\} \in L_0[q_n]}}
\frac{\sum_{i=1}^{n} E_{\boldsymbol{\theta}_{0i}} \lVert\widehat{\boldsymbol{\theta}}^R_i - \boldsymbol{\theta}_{0i}\rVert_{\Sigma}^2} 
{2 q_n \log(n/q_n)} \leq 1.
\end{equation}
In particular, since the minimax error \eqref{minimaxrate} provides a lower bound, one gets the result 
\begin{equation}\label{thm3.second}
\lim\limits_{n \rightarrow \infty} { \sup_{\{\boldsymbol{\theta}_{0i}\} \in L_0[q_n]}}
\frac{\sum_{i=1}^{n} E_{\boldsymbol{\theta}_{0i}} \lVert\widehat{\boldsymbol{\theta}}^R_i - \boldsymbol{\theta}_{0i}\rVert_{\Sigma}^2} 
{2 q_n \log(n/q_n)} = 1.
\end{equation}
\end{theorem}

The following theorem provides results on the rates of posterior contraction for this prior around both the Bayes estimators and the true means.
By \eqref{contraction_est}, the posterior distributions contracts around the Bayes estimator at least as fast as at the minimax rate.
However, by \eqref{contraction_truth}, the rate of posterior contraction around the true means would be slower than the minimax rate.

\begin{theorem}\label{GL_contraction}
    Under the assumptions of Theorem \ref{thm:GL}, we have
    \begin{equation} \label{contraction_est}
         \lim\limits_{n \rightarrow \infty} \sup_{\{\boldsymbol{\theta}_{0i}\} \in L_0[q_n]} E_{\{\boldsymbol{\theta}_{0i}\}} \Pi(\sum_{i=1}^{n} \lVert \boldsymbol{\theta}_i - \widehat{\boldsymbol{\theta}}^R_i  \rVert^2 > q_n \log(\frac{n}{q_n}) \mid \{\boldsymbol{X}_i\}) = 0,
    \end{equation}
    and
    \begin{equation} \label{contraction_truth}
         \limsup\limits_{n \rightarrow \infty} \sup_{\{\boldsymbol{\theta}_{0i}\} \in L_0[q_n]} E_{\{\boldsymbol{\theta}_{0i}\}} \Pi(\sum_{i=1}^{n} \lVert \boldsymbol{\theta}_i - {\boldsymbol{\theta}}_{0i}  \rVert^2 > M_n q_n \log(\frac{n}{q_n}) \mid \{\boldsymbol{X}_i\}) = 0,
    \end{equation}
    for any $\{M_n\}$ such that $\lim_{n\rightarrow\infty} M_n = \infty$.
\end{theorem}

Next we extend the work of \citet{Bhattacharya2015} in the present multivariate framework.
We consider the following prior in which, while (i) remains the same in our earlier formulation, we replace (ii) and (iii) respectively by
\begin{itemize}
    \item[(ii')] $\boldsymbol{\theta}_{i} \vert \lambda_i^2, \tau_i \overset{ind}{\sim} N_k(\boldsymbol{0}, \lambda_i^2 \tau_i\boldsymbol{\Sigma}) $, $i = 1, \dots,n $;
    \item[(iii')] $\lambda_i^2$ and $ \tau_i$ are mutually independent. Also, $\lambda_i^2$'s are independent with $\pi (\lambda_i^2) \propto \exp(-\lambda_i^2/2)$, $i = 1, \dots,n $, while $ \tau_i$'s are also independet with $\pi (\tau_i) \propto \exp(-c_n/(2\tau_i)) \tau_i^{-d-1}$, where $c_n \rightarrow 0$ and will be chosen later and $0<d<1$.
\end{itemize}
As noted by \citet{Bhattacharya2015} as well, the Dirichlet-Laplace priors can be rewritten in the above formulation, except for that the authors put a Gamma prior on $\tau_i$ while we put an Inverse-Gamma prior on it. Due to this discrepancy, we refer the prior defined by (ii') and (iii') as the Exponential-Inverse-Gamma prior.

It is worth mentioning that writing $u_i = \lambda_i^2 \tau_i$, one gets $\pi(u) \propto (u_i + c_n)^{-d-1}$, and one can directly use the $u_i$ for inferential purposes.
Further, this particular formulation is a special case of \citet{vanderPas2016}, who has a very general result concerning asymptotic minimaxity of univariate normal means.
However, it seems more convenient to work with separate priors for $\lambda_i^2$ and $ \tau_i$, and the explicit nature of these priors makes the calculation smooth.
As an aside, $u_i/c_n$ has a beta prime prior with $a=1$ and $b=d$, and this is the prior considered in
\citet{Armagan2011} and \citet{griffin2017hierarchical}.

With the above formulation, the estimators of the  ${\boldsymbol{\theta}}_i$ are denoted as $\widehat{\boldsymbol{\theta}}^{EIG}_i$.
We will prove Theorems \ref{thm:mEIG} in this setup. 
It shows that the prior attains the exact asymptotic minimax error as well.

\begin{theorem}\label{thm:mEIG}
Assume that the true means are sparse in the nearly-black sense.
Under the regularity assumption (II), using the Exponential-Inverse-Gamma prior, i.e., a model satisfying (i), (ii'), and (iii') above,
if $c_n = \left(q_n/n\right)^{\frac{1+\epsilon_n}{d}}$, where $\epsilon_n = 1 / \log(\log(n/q_n))$, then, 
\begin{equation}\label{mEIGThm1}
\lim\limits_{n \rightarrow \infty} { \sup_{\{\boldsymbol{\theta}_{0i}\} \in L_0[q_n]}}{\sum\limits_{i: \boldsymbol{\theta}_{0i} = \boldsymbol{0}} E_0 \lVert \widehat{\boldsymbol{\theta}}^{EIG}_i \rVert^2} /
{\left( q_n \log\left(\frac{n}{q_n}\right) \right)} = 0.
\end{equation}
and
\begin{equation}\label{mEIGThm2}
\limsup\limits_{n \rightarrow \infty} { \sup_{\{\boldsymbol{\theta}_{0i}\} \in L_0[q_n]}}
\frac{\sum\limits_{i: \boldsymbol{\theta}_{0i} \neq \boldsymbol{0}} E_{\boldsymbol{\theta}_{0i}} \lVert\widehat{\boldsymbol{\theta}}^{EIG}_i - \boldsymbol{\theta}_{0i}\rVert^2} 
{2 q_n \log(n/q_n)} \leq 1.
\end{equation}
Consequently,
\begin{equation}\label{mEIGThm3}
\lim\limits_{n \rightarrow \infty} { \sup_{\{\boldsymbol{\theta}_{0i}\} \in L_0[q_n]}}
\frac{\sum_{i=1}^{n} E_{\boldsymbol{\theta}_{0i}} \lVert\widehat{\boldsymbol{\theta}}^{EIG}_i - \boldsymbol{\theta}_{0i}\rVert^2} 
{2 q_n \log(n/q_n)} = 1.
\end{equation}
\end{theorem}

We also have the following results regarding the posterior contraction rate around the Bayes estimator and the true means.
The same contraction rates are observed as using the tuning parameter model.
\begin{theorem}\label{mEIG_contraction}
    Under the assumptions of Theorem \ref{thm:mEIG}, we have
    \begin{equation} \label{contraction_est_mEIG}
         \lim\limits_{n \rightarrow \infty} \sup_{\{\boldsymbol{\theta}_{0i}\} \in l_0[q_n]} E_{\{\boldsymbol{\theta}_{0i}\}} \Pi(\sum_{i=1}^{n} \lVert \boldsymbol{\theta}_i - \widehat{\boldsymbol{\theta}}^{EIG}_i  \rVert^2 > q_n \log(\frac{n}{q_n}) \mid \{\boldsymbol{X}_i\}) = 0,
    \end{equation}
    and
    \begin{equation} \label{contraction_truth_mEIG}
         \limsup\limits_{n \rightarrow \infty} \sup_{\{\boldsymbol{\theta}_{0i}\} \in l_0[q_n]} E_{\{\boldsymbol{\theta}_{0i}\}} \Pi(\sum_{i=1}^{n} \lVert \boldsymbol{\theta}_i - {\boldsymbol{\theta}}_{0i}  \rVert^2 > M_n q_n \log(\frac{n}{q_n}) \mid \{\boldsymbol{X}_i\}) = 0,
    \end{equation}
    for any $\{M_n\}$ such that $\lim_{n\rightarrow\infty} M_n = \infty$.
\end{theorem}

\section{Credible Sets of Multivariate Normal Means}\label{sec:sets}
In this section, we first study coverage probabilities of credible sets constructed under global-local priors defined by (ii) and (iii).
The global parameter is treated as a tuning parameter.
We consider credible sets of the form:
\begin{equation}\label{credible.sets}
    \widehat{C}^{R}_{n,i} = \{\boldsymbol{\theta}_i: \lVert \boldsymbol{\theta}_i - \widehat{\boldsymbol{\theta}}^R_{i} \rVert_{\Sigma}^2 \leq L \widehat{r}_{n,i}^{a/(1+\rho)}(\alpha, \tau_n) \},
\end{equation}
for some $\rho (> 0)$ to be chosen later and $\widehat{r}_{n,i}(\alpha, \tau_n)$ is determined from
\begin{equation*}
    \Pi( \lVert \boldsymbol{\theta}_i - \widehat{\boldsymbol{\theta}}^R_{i} \rVert_{\Sigma}^2 \leq \widehat{r}_{n,i}(\alpha, \tau_n) \mid \boldsymbol{X}_i) = 1 - \alpha.
\end{equation*}
In the following, we will omit the subscript $n$ in $\widehat{C}^{R}_{n,i}$ and $\widehat{r}_{n,i}$ for notational simplicity.

Following \citet{vanDerPas2017}, we view the true mean vectors as in three categories:
\begin{equation*}
    \begin{split}
        \mathcal{S} &:= \{\boldsymbol{\theta}_{0i}: \lVert \boldsymbol{\theta}_{0i} \rVert_{\Sigma}^2 \leq K_S \tau_n\} \\
        \mathcal{M} &:= \{\boldsymbol{\theta}_{0i}: f_{\tau_n} \tau_n \leq \lVert \boldsymbol{\theta}_{0i} \rVert_{\Sigma}^2 \leq K_M \log \frac{1}{\tau_n}\}\\
        \mathcal{L} &:= \{\boldsymbol{\theta}_{0i}: \lVert \boldsymbol{\theta}_{0i} \rVert_{\Sigma}^2 \geq K_L \log \frac{1}{\tau_n}\}
    \end{split}
\end{equation*}
for some positive constants $K_S$, $K_M$ and $K_L$, and some $f_{\tau_n}$ that goes to infinity as $\tau_n$ goes to zero.
We will show that, the proposed credible sets will cover the true means in either $\mathcal{S}$ or $\mathcal{L}$ with a desired probability, while the true means in $\mathcal{M}$ will not be covered with probability tending to one.
The results are summaried in the following theorem.

\begin{theorem}\label{CSThm}
    Consider the global-local prior with a tuning parameter $\tau_n$, i.e., a model satisfying (i), (ii) and (iii), with $a < 1$, under the regularity assumptions (I) and (II).
    Suppose that $K_S > 0$, $K_M < 2a$ and $K_L > 2a$, and that $f_{\tau_n} \rightarrow \infty$ and $f_{\tau_n} \tau_n \rightarrow 0$ as $\tau_n \rightarrow 0$. 
    Then, given $\alpha$, for the credible sets of form \eqref{credible.sets} with $L > \chi^2_{k,\alpha}(\chi^2_{k,\beta})^{-a/(1+\rho)}$ for some fixed $\beta > \alpha$ and $\rho > 0$, 
    \begin{eqnarray}
        &P_{\boldsymbol{\theta}_{0i}} (\boldsymbol{\theta}_{0i} \in \widehat{C}^{R}_{i}) \geq 1-\alpha, \text{ if } \boldsymbol{\theta}_{0i} \in \mathcal{S}, \label{thm.S} \\
            &P_{\boldsymbol{\theta}_{0i}} (\boldsymbol{\theta}_{0i} \notin \widehat{C}^{R}_{i}) \rightarrow 1, \text{ if } \boldsymbol{\theta}_{0i} \in \mathcal{M}, \label{thm.M} \\
            &P_{\boldsymbol{\theta}_{0i}} (\boldsymbol{\theta}_{0i} \in \widehat{C}^{R}_{i}) \geq 1-\alpha, \text{ if } \boldsymbol{\theta}_{0i} \in \mathcal{L}. \label{thm.L}
    \end{eqnarray}
    as $\tau_n \rightarrow 0$.
\end{theorem}
\begin{remark}
    From the proof of the theorem, the conclusions for $\boldsymbol{\theta}_{0i}$ in either $\mathcal{S}$ or $\mathcal{M}$ do not rely on any specific choice of $L$, while only the conclusion for $\mathcal{L}$ has the requirement on $L$.
    To make the credible sets as narrow as possible, noticing that $L > \chi^2_{k,\alpha}(\chi^2_{k,\beta})^{-a/(1+\rho)}$,
    we should choose $\beta$ to be as close to $\alpha$ as possible.
    As for $\rho$, noticing that 
    $$L \widehat{r}_i^{a/(1+\rho)}(\alpha, \tau_n) > \chi^2_{k,\alpha}(\widehat{r}_i(\alpha, \tau_n)/\chi^2_{k,\beta})^{a/(1+\rho)},$$ 
    the choice should depend on $\widehat{r}_i(\alpha, \tau_n)/\chi^2_{k,\beta}$.
    For instance, when $\widehat{r}_i(\alpha, \tau_n)/\chi^2_{k,\beta} > 1$, 
    we can choose as large $\rho$ as possible, so that $L \widehat{r}_i^{a/(1+\rho)}(\alpha, \tau_n)$ would essentially become $\chi^2_{k,\alpha}$. On the other hand, when $\widehat{r}_i(\alpha, \tau_n)/\chi^2_{k,\beta} < 1$, it would be more preferable to choose $\rho$ closer to $0$.
    This observation also motivates an individualized choice of $L_i$ instead of a common $L$ among all the subjects, so that each credible set can be narrowed as much as possible while maintaining the theoretical coverage probability.
\end{remark}

Assuming sparsity in the nearly-black sense, most true means would be in the set $\mathcal{S}$. This fact immediately leads to a high overall coverage probability, i.e., the following corollary.
\begin{corollary}
    Under the setup of Theorem \ref{CSThm}, further assume that the true means $\boldsymbol{\theta}_{0i}$ are sparse in the nearly-black sense. 
    Then, for almost all $i=1,\dots,n$, as $\tau_n \rightarrow 0$, 
    \begin{equation*}
        P_{\boldsymbol{\theta}_{0i}} (\boldsymbol{\theta}_{0i} \in \widehat{C}^{R}_{i}) \geq 1-\alpha.
    \end{equation*}
\end{corollary}

Next, we study coverage probabilities of credible sets constructed under the Exponential-Inverse-Gamma priors defined by (ii') and (iii').
We consider credible sets of the same form as in the previous setup:
\begin{equation}\label{credible.sets.EIG}
    \widehat{C}^{EIG}_{i} = \{\boldsymbol{\theta}_i: \lVert \boldsymbol{\theta}_i - \widehat{\boldsymbol{\theta}}^{EIG}_{i} \rVert_{\Sigma}^2 \leq L \widehat{r}_{i}^{d/(1+\rho)}(\alpha, c_n) \},
\end{equation}
for some $\rho (> 0)$ to be chosen later and $\widehat{r}_{i}(\alpha, c_n)$ is determined from
\begin{equation*}
    \Pi( \lVert \boldsymbol{\theta}_i - \widehat{\boldsymbol{\theta}}^{EIG}_{i} \rVert_{\Sigma}^2 \leq \widehat{r}_{i}(\alpha, c_n) \mid \boldsymbol{X}_i) = 1 - \alpha.
\end{equation*}
Here, we divide the true mean vectors as in the following three categories:
\begin{equation*}
    \begin{split}
        \mathcal{S}' &:= \{\boldsymbol{\theta}_{0i}: \lVert \boldsymbol{\theta}_{0i} \rVert_{\Sigma}^2 \leq K'_S c_n\} \\
        \mathcal{M}' &:= \{\boldsymbol{\theta}_{0i}: f'_{c_n} c_n \leq \lVert \boldsymbol{\theta}_{0i} \rVert_{\Sigma}^2 \leq K'_M \log \frac{1}{c_n}\}\\
        \mathcal{L}' &:= \{\boldsymbol{\theta}_{0i}: \lVert \boldsymbol{\theta}_{0i} \rVert_{\Sigma}^2 \geq K'_L \log \frac{1}{c_n}\}
    \end{split}
\end{equation*}
for some positive constants $K'_S$, $K'_M$ and $K'_L$, and some $f'_{c_n}$ that goes to infinity as $c_n$ goes to zero.
And similar results regarding the coverage probabilities are observed under this prior.

\begin{theorem}\label{CSThm.EIG}
    Consider the Exponential-Inverse-Gamma prior, i.e., a model satisfying (i), (ii') and (iii'), under the regularity assumption (II).
    Suppose that $K'_S > 0$, $K'_M < 2d$ and $K'_L > 2d$, and that $f'_{c_n} \rightarrow \infty$ and $f'_{c_n} c_n \rightarrow 0$ as $c_n \rightarrow 0$. 
    Then, given $\alpha$, for the credible sets of form \eqref{credible.sets.EIG} with $L > \chi^2_{k,\alpha}(\chi^2_{k,\beta})^{-d/(1+\rho)}$ for some fixed $\beta > \alpha$ and $\rho > 0$, for $\alpha < 1/2$
    \begin{eqnarray}
        &P_{\boldsymbol{\theta}_{0i}} (\boldsymbol{\theta}_{0i} \in \widehat{C}^{EIG}_{i}) \geq 1-\alpha, \text{ if } \boldsymbol{\theta}_{0i} \in \mathcal{S'}, \label{thm.S.EIG} \\
            &P_{\boldsymbol{\theta}_{0i}} (\boldsymbol{\theta}_{0i} \notin \widehat{C}^{EIG}_{i}) \rightarrow 1, \text{ if } \boldsymbol{\theta}_{0i} \in \mathcal{M'}, \label{thm.M.EIG} \\
            &P_{\boldsymbol{\theta}_{0i}} (\boldsymbol{\theta}_{0i} \in \widehat{C}^{EIG}_{i}) \geq 1-\alpha, \text{ if } \boldsymbol{\theta}_{0i} \in \mathcal{L'}. \label{thm.L.EIG}
    \end{eqnarray}
    as $\tau_n \rightarrow 0$.
\end{theorem}
The following corollary is also immediate due to the nearly-black sparsity.
\begin{corollary}
    Under the setup of Theorem \ref{CSThm.EIG}, further assume that the true means $\boldsymbol{\theta}_{0i}$ are sparse in the nearly-black sense. 
    Then, for almost all $i=1,\dots,n$,  as $c_n \rightarrow 0$, 
    \begin{equation*}
        P_{\boldsymbol{\theta}_{0i}} (\boldsymbol{\theta}_{0i} \in \widehat{C}^{EIG}_{i}) \geq 1-\alpha.
    \end{equation*}
\end{corollary}

\section{Final Remarks}
\label{sec:remark}
The paper addresses asymptotic estimation of multivariate normal means under global-local priors.
We first find the asymptotic minimax error in the multivariate setup.
Then, the asymptotic minimax error is obtained by treating the global parameter as a tuning parameter.
The same result is obtained under Dirichlet-Laplace priors. 
Also, credible sets are obtained under global-local priors extending the idea of \citet{vanDerPas2017} in the multivariate case. 

\begin{appendices}

\section{Lemmas}
\label{sec:appendixa}
\subsection{Lemmas for the Multivariate Tuning Parameter Model}

Regarding the multivariate tuning parameter model, we can establish the following lemmas.
Under the model (i) - (iii), the posterior density of $\kappa_i$ is given by
\begin{equation*}
    \begin{split}
        \Pi(\kappa_i\vert X_i) &\propto \kappa^{k/2+a-1} (1-\kappa)^{-a-1} L\left(\frac{1-\kappa}{\kappa \tau_n}\right) 
    \exp{\left(-\frac{\kappa}{2}  \boldsymbol{X}_i^T \boldsymbol{\Sigma}^{-1} \boldsymbol{X}_i  \right)}.
    \end{split}
\end{equation*}
In the following, we will use $K (> 0)$ to denote a generic constant.
\begin{lemma}\label{TUNElemma1}
Under the multivariate global-local prior model treating the global parameter as a tuning parameter and under the regularity assumption (I), 
\begin{equation}
    E(1-\kappa_i\vert \boldsymbol{X}_i) \le K \tau_n^a \exp \left(\frac{ \boldsymbol{X}_i^T \boldsymbol{\Sigma}^{-1} \boldsymbol{X}_i }{2}\right).
\end{equation}
\end{lemma}
\begin{proof}
For an arbitrary constant $\xi>0$,
\begin{equation}\label{3}
\begin{split}
    &E\left(1-\kappa_i  \mid  \boldsymbol{X}_i\right) \\
    = &\tau_n E\left[\lambda_i^2 \left(1 + \lambda_i^2 \tau_n\right)^{-1}  \mid  \boldsymbol{X}_i\right] \\
    = &\tau_n \frac{\int_0^\infty \left(1+\lambda^2 \tau_n\right)^{-\left(\frac{k}{2}+1\right)} \left(\lambda^2\right)^{-a} L\left(\lambda^2\right) 
    \exp \left(-\frac{ \boldsymbol{X}_i^T \boldsymbol{\Sigma}^{-1} \boldsymbol{X}_i }{2\left(1+\lambda^2 \tau_n \right)}\right)d\lambda^2}
    {\int_0^\infty (1+\lambda^2 \tau_n)^{-\frac{k}{2}} (\lambda^2)^{-a-1} L(\lambda^2) 
    \exp \left(-\frac{ \boldsymbol{X}_i^T \boldsymbol{\Sigma}^{-1} \boldsymbol{X}_i }{2(1+\lambda^2 \tau_n)}\right)d\lambda^2} \\
    \le &\tau_n exp(\frac{ \boldsymbol{X}_i^T \boldsymbol{\Sigma}^{-1} \boldsymbol{X}_i }{2}) I_1 / \left(I_2 L(\xi)\right)
\end{split}
\end{equation}
where 
\begin{equation}
    I_1  = \int_0^\infty \left(1+\lambda^2 \tau_n\right)^{-\left(\frac{k}{2}+1\right)} \left(\lambda^2\right)^{-a} L\left(\lambda^2\right) 
    d\lambda^2
\end{equation}
and 
\begin{equation}
    I_2 = \int_{\xi}^\infty \left(1+\lambda^2 \tau_n\right)^{-\frac{k}{2}} \left(\lambda^2\right)^{-a-1}
    d\lambda^2
\end{equation}

But 
\begin{equation}\label{4}
    I_2 \ge 
    \int_{\xi}^\infty \left(1+\lambda^2\right)^{-\frac{k}{2}} \left(\lambda^2\right)^{-a-1} d\lambda^2
    = K.
\end{equation}

When $0 < a < 1$,
\begin{equation}
\begin{split}
    I_1 \le M \tau_n^{a-1} \int_0^\infty \left(1+\lambda^2 \tau_n\right)^{-\left(\frac{k}{2}+1\right)} \left(\tau_n \lambda^2\right)^{-a} d\tau_n \lambda^2 
    = K \tau_n^{a-1}.
\end{split}
\end{equation}
Combining this with \eqref{3}-\eqref{4},  $E\left(1-\kappa_i  \mid  \boldsymbol{X}_i\right) \le K \tau_n^a \exp \left(\frac{ \boldsymbol{X}_i^T \boldsymbol{\Sigma}^{-1} \boldsymbol{X}_i }{2}\right)$, when  $0 < a < 1$.

When $a=1$,
\begin{equation}
\begin{split}
    I_1 &= \int_0^\infty \left(1+\lambda^2 \tau_n\right)^{-\left(\frac{k}{2}+1\right)} \left(\tau_n \lambda^2\right)^{-1} L(\lambda^2) d\tau_n \lambda^2  \\
    &= \int_0^\infty \left(1+x \right)^{-\left(\frac{k}{2}+1\right)} x^{-1} L(x / \tau_n) dx \\
    &= \left(\int_0^1 + \int_1^\infty\right) \left(1+x \right)^{-\left(\frac{k}{2}+1\right)} x^{-1} L(x / \tau_n) dx,
\end{split}
\end{equation}
in which, under regularity assumption (I),
\begin{equation}
\begin{split}
    \int_0^1 \left(1+x \right)^{-\left(\frac{k}{2}+1\right)} x^{-1} L(x / \tau_n) dx &\leq \int_0^1 x^{-2}  L(x / \tau_n) L(x) / L(x) dx \\
    &\leq \frac{M}{m} \int_0^1 x^{-2} L(x) dx = K,
\end{split}
\end{equation}
and
\begin{equation}
\begin{split}
    \int_1^\infty \left(1+x \right)^{-\left(\frac{k}{2}+1\right)} x^{-1} L(x / \tau_n) dx \leq M \int_1^\infty \left(1+x \right)^{-1} x^{-1}dx = K.
\end{split}
\end{equation}
Hence, when $a = 1$, $E\left(1-\kappa_i  \mid  \boldsymbol{X}_i\right) \le K \tau_n \exp \left(\frac{ \boldsymbol{X}_i^T \boldsymbol{\Sigma}^{-1} \boldsymbol{X}_i }{2}\right)$.

Lastly, when $a > 1$, we write
\begin{equation}
\begin{split}
    I_1 = \tau_n^{a-1} \left(\int_0^1 + \int_1^\infty\right) \left(1+x \right)^{-\left(\frac{k}{2}+1\right)} x^{-a} L(x / \tau_n) dx,
\end{split}
\end{equation}
Then, similar techeniques lead to that $E\left(1-\kappa_i  \mid  \boldsymbol{X}_i\right) \le K \tau_n^a \exp \left(\frac{ \boldsymbol{X}_i^T \boldsymbol{\Sigma}^{-1} \boldsymbol{X}_i }{2}\right)$.
Therefore, we know the inequality holds for all $a > 0$.
\end{proof}

\begin{lemma}\label{TUNElemma2}
Under the multivariate global-local prior model with treating the global parameter as a tuning parameter and under the regularity assumption (I), for arbitrary constants $0 < \xi < 1$ and $0 < \delta < 1$,
\begin{equation}
    E \left(\kappa_i \mathbbm{1}_{\left[ \kappa_i > \xi \right]} \mid  \boldsymbol{X}_i\right) 
    \leq K \tau_n^{-a} \exp{\left[ -\frac{\xi(1-\delta)}{2}  \boldsymbol{X}_i^T \boldsymbol{\Sigma}^{-1} \boldsymbol{X}_i  \right]}.
\end{equation}
\end{lemma}
\begin{proof}
For an arbitrary constant $0 < \xi < 1$,
\begin{equation}\label{13}
\begin{split}
    &E \left(\kappa_i \mathbbm{1}_{\left[ \kappa_i > \xi \right]} \mid  \boldsymbol{X}_i\right) \\
    =& \frac{\int_{\xi}^1 \kappa^{k/2+a} (1-\kappa)^{-a-1} L\left(\frac{1-\kappa}{\kappa \tau_n}\right) 
    \exp{\left(-\frac{\kappa}{2}  \boldsymbol{X}_i^T \boldsymbol{\Sigma}^{-1} \boldsymbol{X}_i  \right)} d\kappa}
    {\int_0^{1} \kappa^{k/2+a-1} (1-\kappa)^{-a-1} L\left(\frac{1-\kappa}{\kappa \tau_n}\right) 
    \exp{\left(-\frac{\kappa}{2}  \boldsymbol{X}_i^T \boldsymbol{\Sigma}^{-1} \boldsymbol{X}_i  \right)} d\kappa} \\
    \leq& \frac{\int_{\xi}^1 \kappa^{k/2+a} (1-\kappa)^{-a-1} L\left(\frac{1-\kappa}{\kappa \tau_n}\right) 
    \exp{\left(-\frac{\kappa}{2}  \boldsymbol{X}_i^T \boldsymbol{\Sigma}^{-1} \boldsymbol{X}_i  \right)} d\kappa}
    {\int_0^{\xi\delta} \kappa^{k/2+a-1} (1-\kappa)^{-a-1} L\left(\frac{1-\kappa}{\kappa \tau_n}\right) 
    \exp{\left(-\frac{\kappa}{2}  \boldsymbol{X}_i^T \boldsymbol{\Sigma}^{-1} \boldsymbol{X}_i  \right)} d\kappa} \\
    \leq& \exp{\left[-\frac{\xi(1-\delta)}{2} \boldsymbol{X}_i^T \boldsymbol{\Sigma}^{-1} \boldsymbol{X}_i \right]}
    \frac{ \int_{\xi}^1 \kappa^{k/2+a} (1-\kappa)^{-a-1} L\left(\frac{1-\kappa}{\kappa \tau_n}\right) d\kappa}
    {\int_0^{\xi\delta} \kappa^{k/2+a-1} (1-\kappa)^{-a-1} L\left(\frac{1-\kappa}{\kappa \tau_n}\right) d\kappa}.
\end{split}
\end{equation}

Now observe that
\begin{equation}\label{14}
\begin{split}
    &\int_{\xi}^1 \kappa^{k/2+a} (1-\kappa)^{-a-1} L\left(\frac{1-\kappa}{\kappa \tau_n}\right) d\kappa \\
    =&\int_{0}^{(1-\xi)/(\xi \tau_n)} (1+\lambda^2 \tau_n)^{-(k/2+a)} (\lambda^2 \tau_n (1+\lambda^2 \tau_n)^{-1})^{-a-1} L(\lambda^2) \frac{\tau_n d \lambda^2}{(1+\lambda^2 \tau_n)^{2}} \\
    =& \tau_n^{-a} \int_{0}^{(1-\xi)/(\xi \tau_n)} (1+\lambda^2 \tau_n)^{-(k/2+1)} (\lambda^2)^{-a-1} L(\lambda^2) d \lambda^2 \\
    \leq& \tau_n^{-a} \int_0^\infty (\lambda^2)^{-a-1} L(\lambda^2) d \lambda^2 = \tau_n^{-a} K^{-1}.
\end{split}
\end{equation}

By assumption (I) and $\tau_n<1$,
\begin{equation}\label{15}
\begin{split}
    &\int_0^{\xi\delta} \kappa^{k/2+a-1} (1-\kappa)^{-a-1} L\left(\frac{1-\kappa}{\kappa \tau_n}\right) d\kappa \\
    \geq&\int_{0}^{\xi\delta} \kappa^{k/2+a-1} L\left(\frac{1-\kappa}{\kappa \tau_n}\right) d\kappa \\
    \geq& L\left(\frac{1-\xi\delta}{\xi\delta}\right) \int_{0}^{\xi\delta} \kappa^{k/2+a-1} d \kappa \\
    =& L\left(\frac{1-\xi\delta}{\xi\delta}\right) \frac{(\xi\delta)^{k/2+a}}{k/2+a}.  
\end{split}
\end{equation}

Combining \eqref{13} - \eqref{15},
\begin{equation}
    E \left(\kappa_i \mathbbm{1}_{\left[ \kappa_i > \xi \right]} \mid  \boldsymbol{X}_i\right) 
    \leq K \tau_n^{-a} \exp{\left[ -\frac{\xi(1-\delta)}{2}  \boldsymbol{X}_i^T \boldsymbol{\Sigma}^{-1} \boldsymbol{X}_i  \right]}.
\end{equation}
\end{proof}

\begin{lemma}\label{TUNElemma3}
Under the multivariate global-local prior model with treating the global parameter as a tuning parameter and under the regularity assumption (I), for an arbitrary constant $0 < \xi < 1$,
\begin{equation}
    E(\kappa_i \mathbbm{1}_{[\kappa_i \le \xi]} \vert \boldsymbol{X}_i) \le K / \boldsymbol{X}_i^T\boldsymbol{\Sigma}^{-1}\boldsymbol{X}_i.
\end{equation}
\end{lemma}
\begin{proof}
For an arbitrary constant $0 < \xi < 1$,
\begin{equation}\label{11}
\begin{split}
    &E \left(\kappa_i \mathbbm{1}_{\left[ \kappa_i \leq \xi \right]} \mid  \boldsymbol{X}_i\right) \\
    =& \frac{\int_0^{\xi} \kappa^{k/2+a} (1-\kappa)^{-a-1} L\left(\frac{1-\kappa}{\kappa \tau_n}\right) 
    \exp{\left(-\frac{\kappa}{2}  \boldsymbol{X}_i^T \boldsymbol{\Sigma}^{-1} \boldsymbol{X}_i  \right)} d\kappa}
    {\int_0^{1} \kappa^{k/2+a-1} (1-\kappa)^{-a-1} L\left(\frac{1-\kappa}{\kappa \tau_n}\right) 
    \exp{\left(-\frac{\kappa}{2}  \boldsymbol{X}_i^T \boldsymbol{\Sigma}^{-1} \boldsymbol{X}_i  \right)} d\kappa} \\
    \leq& \frac{(1-\xi)^{-a-1}M}{L\left(\frac{1-\xi}{\xi \tau_n}\right)}
    \frac{\int_0^{\xi} \kappa^{k/2+a} \exp{\left(-\frac{\kappa}{2}  \boldsymbol{X}_i^T \boldsymbol{\Sigma}^{-1} \boldsymbol{X}_i  \right)} d\kappa}
    {\int_0^{\xi} \kappa^{k/2+a-1} \exp{\left(-\frac{\kappa}{2}  \boldsymbol{X}_i^T \boldsymbol{\Sigma}^{-1} \boldsymbol{X}_i  \right)} d\kappa} \\
    =& \frac{(1-\xi)^{-a-1}M}{L\left(\frac{1-\xi}{\xi \tau_n}\right)}
    \frac{\int_0^{\xi  \boldsymbol{X}_i^T \boldsymbol{\Sigma}^{-1} \boldsymbol{X}_i /2} \left(\frac{2t}{ \boldsymbol{X}_i^T \boldsymbol{\Sigma}^{-1} \boldsymbol{X}_i }\right)^{k/2+a} \exp{\left(-t\right)} dt}
    {\int_0^{\xi  \boldsymbol{X}_i^T \boldsymbol{\Sigma}^{-1} \boldsymbol{X}_i /2} \left(\frac{2t}{ \boldsymbol{X}_i^T \boldsymbol{\Sigma}^{-1} \boldsymbol{X}_i }\right)^{k/2+a-1} \exp{\left(-t\right)} dt} \\
    =& \frac{(1-\xi)^{-a-1}M}{L\left(\frac{1-\xi}{\xi \tau_n}\right)  \boldsymbol{X}_i^T \boldsymbol{\Sigma}^{-1} \boldsymbol{X}_i } 
    \frac{\int_0^{\xi  \boldsymbol{X}_i^T \boldsymbol{\Sigma}^{-1} \boldsymbol{X}_i /2} t^{k/2+a} \exp{\left(-t\right)} dt}
    {\int_0^{\xi  \boldsymbol{X}_i^T \boldsymbol{\Sigma}^{-1} \boldsymbol{X}_i /2} t^{k/2+a-1} \exp{\left(-t\right)} dt}.
\end{split}
\end{equation}

Integrating the numerator by parts,
\begin{equation}\label{12}
\begin{split}
    &\int_0^{\xi  \boldsymbol{X}_i^T \boldsymbol{\Sigma}^{-1} \boldsymbol{X}_i /2} t^{k/2+a} \exp{\left(-t\right)} dt = (k/2 + a) \times \\
    &\int_0^{\xi  \boldsymbol{X}_i^T \boldsymbol{\Sigma}^{-1} \boldsymbol{X}_i /2} t^{k/2+a-1} \exp{\left(-t\right)} dt - 
    \left( \frac{\xi \boldsymbol{X}_i^T \boldsymbol{\Sigma}^{-1} \boldsymbol{X}_i }{2} \right)^{\frac{k}{2}+a} \exp{\left( -\frac{\xi \boldsymbol{X}_i^T \boldsymbol{\Sigma}^{-1} \boldsymbol{X}_i }{2} \right)}.
\end{split}
\end{equation}

Combining \eqref{11} and \eqref{12},
\begin{equation}\label{ineq.lemma1}
    E \left(\kappa_i \mathbbm{1}_{\left[ \kappa_i \leq \xi \right]} \mid  \boldsymbol{X}_i\right) \leq K/ \boldsymbol{X}_i^T \boldsymbol{\Sigma}^{-1} \boldsymbol{X}_i .
\end{equation}
\end{proof}

\subsection{Lemmas for the Multivariate Exponential-Inverse-Gamma Model}
Under an Exponential-Inverse-Gamma model (i), (ii') and (iii')
the posterior density of $\kappa_i$ is given by
\begin{equation*}
    \begin{split}
        \Pi(\kappa_i\vert \boldsymbol{X}_i) \propto \kappa_i^{d+\frac{k}{2}-1} (1-\kappa_i + \kappa_i c_n)^{-d-1} \exp\left(-\frac{\kappa_i \boldsymbol{X}_i^T\boldsymbol{\Sigma}^{-1}\boldsymbol{X}_i}{2}\right) .
    \end{split}
\end{equation*}

Now, we give three basic inequalities involving the Exponential-Inverse-Gamma prior. 

\begin{lemma}\label{mEIGlemma1}
Under the multivariate Exponential-Inverse-Gamma model, assuming $c_n \rightarrow 0$
as $n \rightarrow \infty$, for large $n$,
\begin{equation}
    E(1-\kappa_i\vert \boldsymbol{X}_i) \le K c_n^{d}  \exp \left(\frac{\boldsymbol{X}_i^T\boldsymbol{\Sigma}^{-1}\boldsymbol{X}_i}{2}\right).
\end{equation}
\end{lemma}
\begin{proof}
Firstly,
\begin{equation*}
    \begin{split}
        E(1-\kappa_i\vert \boldsymbol{X}_i) &= \frac{\int_0^1 (1-\kappa_i) \kappa_i^{d+\frac{k}{2}-1} \exp\left(-\frac{\kappa_i \boldsymbol{X}_i^T\boldsymbol{\Sigma}^{-1}\boldsymbol{X}_i}{2}\right)  (1-\kappa_i + \kappa_i c_n)^{-d-1} d\kappa_i}{\int_0^1 \kappa_i^{d+\frac{k}{2}-1} \exp\left(-\frac{\kappa_i \boldsymbol{X}_i^T\boldsymbol{\Sigma}^{-1}\boldsymbol{X}_i}{2}\right) (1-\kappa_i + \kappa_i c_n)^{-d-1} d\kappa_i} \\
        &\le \exp(\boldsymbol{X}_i^T\boldsymbol{\Sigma}^{-1}\boldsymbol{X}_i/2)  \frac{\int_0^1 (1-\kappa_i)\kappa_i^{d+\frac{k}{2}-1} (1-\kappa_i + \kappa_i c_n)^{-d-1} d\kappa_i}{\int_0^1 \kappa_i^{d+\frac{k}{2}-1}  (1-\kappa_i + \kappa_i c_n)^{-d-1} d\kappa_i}\\
        & := \exp(\boldsymbol{X}_i^T\boldsymbol{\Sigma}^{-1}\boldsymbol{X}_i/2) \frac{N}{D}.
    \end{split}
\end{equation*}
Now,
\begin{equation*}
    \begin{split}
        D &\ge \int_{1-c_n}^1 \kappa_i^{d+\frac{k}{2}-1}   (1-\kappa_i + \kappa_i c_n)^{-d-1} d\kappa_i\\
        & \ge \int_{1-c_n}^1 \kappa_i^{d+\frac{k}{2}-1}   (c_n +  c_n )^{-d-1} d\kappa_i\\
        & = \left(2 c_n\right)^{-(d+1)} \int_{1-c_n}^1 \kappa_i^{d+\frac{k}{2}-1} d \kappa_i\\
        & = \left(2 c_n\right)^{-(d+1)}\frac{ 1-(1-c_n)^{d+\frac{k}{2}}}{d+\frac{k}{2}}\\
        & \ge K c_n^{-(d+1)} \left[ c_n- \frac{d+\frac{k}{2}-1}{2}c_n^2\right]\\
        & = K c_n^{-d}  \left[1- \frac{d+\frac{k}{2}-1}{2}c_n\right]\\
        & \ge K c_n^{-d}.
    \end{split}
\end{equation*}
Next,
\begin{equation*}
\begin{split}
    N &= \int_0^1 (1-\kappa_i)\kappa_i^{d+\frac{k}{2}-1} (1-\kappa_i + \kappa_i c_n)^{-d-1} d\kappa_i \\
    & \leq \int_0^1 (1-\kappa_i)\kappa_i^{d+\frac{k}{2}-1} (1-\kappa_i)^{-d-1} d\kappa_i \\
    & = \int_0^1 \kappa_i^{d+\frac{k}{2}-1} (1-\kappa_i)^{-d} d\kappa_i \\
    & = \text{Beta}(1-d, \frac{k}{2}+d).
\end{split}
\end{equation*}
Therefore, when $n$ is sufficiently large,
\begin{equation*}
     E(1-\kappa_i \vert \boldsymbol{X}_i) \le  K c_n^{d} \exp(\boldsymbol{X}_i^T\boldsymbol{\Sigma}^{-1}\boldsymbol{X}_i/2) .
\end{equation*}
\end{proof}

\begin{lemma}\label{mEIGlemma2}
Under the multivariate Exponential-Inverse-Gamma model, for arbitrary constants $0 < \xi < 1$ and $0 < \delta < 1$,
\begin{equation}
    E(\kappa_i \mathbbm{1}_{[\kappa_i > \xi]} \vert \boldsymbol{X}_i) \le K c_n^{-d} \exp\left(-\frac{1}{2}\xi(1-\delta)\boldsymbol{X}_i^T\boldsymbol{\Sigma}^{-1}\boldsymbol{X}_i\right).
\end{equation}
\end{lemma}
\begin{proof}
For arbitrary constants $0 < \xi < 1$ and $0 < \delta < 1$,
\begin{equation*}
    \begin{split}
         & E(\kappa_i \mathbbm{1}_{ [\kappa_i > \xi] }\vert \boldsymbol{X}_i)\\
        = &\frac{\int_{\xi}^1 \kappa_i \kappa_i^{d+\frac{k}{2}-1} \exp\left(-\frac{\kappa_i \boldsymbol{X}_i^T\boldsymbol{\Sigma}^{-1}\boldsymbol{X}_i}{2}\right)   (1-\kappa_i + \kappa_i c_n)^{-d-1} d\kappa_i}{\int_0^1 \kappa_i^{d+\frac{k}{2}-1} \exp\left(-\frac{\kappa_i \boldsymbol{X}_i^T\boldsymbol{\Sigma}^{-1}\boldsymbol{X}_i}{2}\right)   (1-\kappa_i + \kappa_i c_n)^{-d-1} d\kappa_i} \\
        \le & \frac{\int_{\xi}^1  \kappa_i^{d+\frac{k}{2}} \exp\left(-\frac{\kappa_i \boldsymbol{X}_i^T\boldsymbol{\Sigma}^{-1}\boldsymbol{X}_i}{2}\right)   (1-\kappa_i + \kappa_i c_n)^{-d-1} d\kappa_i}{\int_0^{\xi\delta} \kappa_i^{d+\frac{k}{2}-1} \exp\left(-\frac{\kappa_i \boldsymbol{X}_i^T\boldsymbol{\Sigma}^{-1}\boldsymbol{X}_i}{2}\right)  (1-\kappa_i + \kappa_i c_n)^{-d-1} d\kappa_i}\\
        \le & \exp\left(-\frac{1}{2} \xi (1-\delta) \boldsymbol{X}_i^T\boldsymbol{\Sigma}^{-1}\boldsymbol{X}_i \right) \frac{\int_{\xi}^1  \kappa_i^{d+\frac{k}{2}}   (1-\kappa_i + \kappa_i c_n)^{-d-1} d\kappa_i}{\int_0^{\xi\delta} \kappa_i^{d+\frac{k}{2}-1}    (1-\kappa_i + \kappa_i c_n)^{-d-1} d\kappa_i}\\
        \le & \exp\left(-\frac{1}{2} \xi (1-\delta) \boldsymbol{X}_i^T\boldsymbol{\Sigma}^{-1}\boldsymbol{X}_i \right) \frac{\int_{\xi}^1  \kappa_i^{d+\frac{k}{2}}   (1-\kappa_i + \kappa_i c_n)^{-d-1} d\kappa_i}{\int_0^{\xi\delta} \kappa_i^{d+\frac{k}{2}-1}    (1 + c_n)^{-d-1} d\kappa_i}\\
        \le & \exp\left(-\frac{1}{2} \xi (1-\delta) \boldsymbol{X}_i^T\boldsymbol{\Sigma}^{-1}\boldsymbol{X}_i \right) (1+c_n)^{d+1} \frac{d+k/2}{(\xi\delta)^{d+k/2}}\int_{\xi}^1   (1-\kappa_i + \kappa_i c_n)^{-d-1} d\kappa_i\\
        \le & K \exp\left(-\frac{1}{2} \xi (1-\delta) \boldsymbol{X}_i^T\boldsymbol{\Sigma}^{-1}\boldsymbol{X}_i \right) \left[ \frac{(1-\kappa_i + \xi c_n)^{-d}}{d}\right]^1_\xi\\
        \le &  K \exp\left(-\frac{1}{2} \xi (1-\delta) \boldsymbol{X}_i^T\boldsymbol{\Sigma}^{-1}\boldsymbol{X}_i \right)c_n^{-d}.
    \end{split}
\end{equation*}
\end{proof}

\begin{lemma}\label{mEIGlemma3}
Under the multivariate Exponential-Inverse-Gamma model, for an arbitrary constant $0 < \xi < 1$,
\begin{equation}
    E(\kappa_i \mathbbm{1}_{[\kappa_i \le \xi]} \vert \boldsymbol{X}_i) \le K / \boldsymbol{X}_i^T\boldsymbol{\Sigma}^{-1}\boldsymbol{X}_i.
\end{equation}
\end{lemma}
\begin{proof}
For an arbitrary constant $0 < \xi < 1$,
\begin{equation*}
    \begin{split}
         & E(\kappa_i \mathbbm{1}_{ [\kappa_i \le \xi] }\vert \boldsymbol{X}_i)\\
        = &\frac{\int_0^{\xi} \kappa_i^{d+\frac{k}{2}} \exp\left(-\frac{\kappa_i \boldsymbol{X}_i^T\boldsymbol{\Sigma}^{-1}\boldsymbol{X}_i}{2}\right)   (1-\kappa_i + \kappa_i c_n)^{-d-1} d\kappa_i}{\int_0^1 \kappa_i^{d+\frac{k}{2}-1} \exp\left(-\frac{\kappa_i \boldsymbol{X}_i^T\boldsymbol{\Sigma}^{-1}\boldsymbol{X}_i}{2}\right)   (1-\kappa_i + \kappa_i c_n)^{-d-1} d\kappa_i} \\
        \le &\frac{\int_0^{\xi} \kappa_i^{d+\frac{k}{2}} \exp\left(-\frac{\kappa_i \boldsymbol{X}_i^T\boldsymbol{\Sigma}^{-1}\boldsymbol{X}_i}{2}\right)   (1-\kappa_i + \kappa_i c_n)^{-d-1} d\kappa_i}{\int_0^\xi \kappa_i^{d+\frac{k}{2}-1} \exp\left(-\frac{\kappa_i \boldsymbol{X}_i^T\boldsymbol{\Sigma}^{-1}\boldsymbol{X}_i}{2}\right)   (1-\kappa_i + \kappa_i c_n)^{-d-1} d\kappa_i} \\
        \le & \frac{(1-\xi)^{-d-1}}{(1+\xi c_n)^{-d-1}}\frac{\int_0^{\xi} \kappa_i^{d+\frac{k}{2}} \exp\left(-\frac{\kappa_i \boldsymbol{X}_i^T\boldsymbol{\Sigma}^{-1}\boldsymbol{X}_i}{2}\right)    d\kappa_i}{\int_0^\xi \kappa_i^{d+\frac{k}{2}-1} \exp\left(-\frac{\kappa_i \boldsymbol{X}_i^T\boldsymbol{\Sigma}^{-1}\boldsymbol{X}_i}{2}\right)d\kappa_i} \\
        \le & \frac{(1-\xi)^{-d-1}}{(1+\xi c_n)^{-d-1} \boldsymbol{X}_i^T\boldsymbol{\Sigma}^{-1}\boldsymbol{X}_i}\frac{\int_0^{\xi \boldsymbol{X}_i^T\boldsymbol{\Sigma}^{-1}\boldsymbol{X}_i/2}  \exp\left(-t\right)   t^{d+\frac{k}{2}} dt}{\int_0^{\xi \boldsymbol{X}_i^T\boldsymbol{\Sigma}^{-1}\boldsymbol{X}_i/2} \exp\left(-t\right) t^{d+\frac{k}{2}-1}dt}. \\
    \end{split}
\end{equation*}
Integrating by parts,
\begin{equation*}
    \begin{split}
        &\int_0^{\xi \boldsymbol{X}_i^T\boldsymbol{\Sigma}^{-1}\boldsymbol{X}_i/2}  \exp\left(-t\right) t^{d+\frac{k}{2}}dt\\
        =& \left[-\exp\left(-t\right)t^{d+\frac{k}{2}}\right]^{t=\xi \boldsymbol{X}_i^T\boldsymbol{\Sigma}^{-1}\boldsymbol{X}_i/2}_{t=0} + (d+\frac{k}{2}) \int_0^{\xi \boldsymbol{X}_i^T\boldsymbol{\Sigma}^{-1}\boldsymbol{X}_i/2}  \exp\left(-t\right)  t^{d+\frac{k}{2}-1}  dt.
    \end{split}
\end{equation*}
Therefore,
\begin{equation*}
     \frac{\int_0^{\xi \boldsymbol{X}_i^T\boldsymbol{\Sigma}^{-1}\boldsymbol{X}_i/2}  \exp\left(-t\right)   ^{d+\frac{k}{2}} dt}{\int_0^{\xi \boldsymbol{X}_i^T\boldsymbol{\Sigma}^{-1}\boldsymbol{X}_i/2} \exp\left(-t\right) t^{d+\frac{k}{2}-1}dt} \le d + \frac{k}{2}.
\end{equation*}
Hence,
\begin{equation*}
 E(\kappa_i \mathbbm{1}_{ [\kappa_i \le \xi] }\vert \boldsymbol{X}_i) \le K/\boldsymbol{X}_i^T\boldsymbol{\Sigma}^{-1}\boldsymbol{X}_i.
\end{equation*}
\end{proof}

\section{Proofs of the Main Theorems}
\label{sec:appendixb}

\subsection{Proofs of Theorems in  Section \ref{sec:multivariate}}
\begin{proof}[Proof of Theorem \ref{thm:minimax.error}]
For each $\boldsymbol{X}_{i}  \overset{ind}{\sim} N_k(\boldsymbol{\theta}_{i}, \boldsymbol{\Sigma})$, define $\boldsymbol{Y}_{i} = (Y_{i1}, \dots, Y_{ik})^\top := \boldsymbol{\Sigma}^{-1/2}\boldsymbol{X}_{i}$ and $\boldsymbol{\lambda}_{i} = (\lambda_{i1}, \dots, \lambda_{ik})^\top := \boldsymbol{\Sigma}^{-1/2}\boldsymbol{\theta}_{i}$.
So, all the components of $\boldsymbol{Y}_{i}$'s are independently normally distributed, i.e.,  ${Y}_{ij}  \overset{ind}{\sim} N({\lambda}_{ij}, 1), i = 1, \dots, n, j=1,\dots,k$.

Given estimators $\widehat{\boldsymbol{\lambda}}_{i}$, if we use $\widehat{\boldsymbol{\theta}}_{i} = \boldsymbol{\Sigma}^{1/2}\widehat{\boldsymbol{\lambda}}_{i}$ to estimate $\boldsymbol{\theta}_{i}$'s, then
$$E \lVert \widehat{\boldsymbol{\theta}}_{i} - \boldsymbol{\theta}_{i} \rVert_{\Sigma}^2 = E \lVert \widehat{\boldsymbol{\lambda}}_{i} - \boldsymbol{\lambda}_{i} \rVert^2 = \sum_{j=1}^{k} E ( \widehat{{\lambda}}_{ij} - {\lambda}_{ij} )^2. $$

Also, since $\{ \lambda_{ij}\} \in l_0[q'_{nk}]$ would imply that $\{ \boldsymbol{\theta}_{i}\} \in L_0[q_n]$, when we let $q'_{nk} = q_n$, we have 
$$\sup_{\{ \lambda_{ij}\} \in l_0[q'_{nk}]} \sum_{i=1}^n \sum_{j=1}^{k} E ( \widehat{{\lambda}}_{ij} - {\lambda}_{ij} )^2 \leq \sup_{\{ \boldsymbol{\theta}_{i}\} \in L_0[q_n]} \sum_{i=1}^n E \lVert \widehat{\boldsymbol{\theta}}_{i} - \boldsymbol{\theta}_{i} \rVert_{\Sigma}^2. $$

Finally, since $\boldsymbol{\Sigma}$ is positive definite, there is a one-to-one correspondence between 
$\{\{\widehat{\boldsymbol{\theta}}_{i}\}_{i=1}^n\}$ and $\{\{\widehat{{\lambda}}_{ij}\}_{i=1,j=1}^{n,\quad k}\}$. So the above inequality still hold, if we further take the infinum over all possible estimators:
$$\inf_{\{\widehat{{\lambda}}_{ij}\}} \sup_{\{ \lambda_{ij}\} \in l_0[q'_{nk}]} \sum_{i=1}^n \sum_{j=1}^{k} E ( \widehat{{\lambda}}_{ij} - {\lambda}_{ij} )^2 \leq \inf_{\{\widehat{\boldsymbol{\theta}}_{i}\}} \sup_{\{ \boldsymbol{\theta}_{i}\} \in L_0[q_n]} \sum_{i=1}^n E \lVert \widehat{\boldsymbol{\theta}}_{i} - \boldsymbol{\theta}_{i} \rVert_{\Sigma}^2. $$

\citet{Donoho1992} provides the result for the left hand side, which is, as $n \rightarrow \infty$,
\begin{equation*}
    \begin{split}
        \inf_{\{\widehat{{\lambda}}_{ij}\}} \sup_{\{ \lambda_{ij}\} \in l_0[q'_{nk}]} \sum_{i=1}^n \sum_{j=1}^{k} E ( \widehat{{\lambda}}_{ij} - {\lambda}_{ij} )^2 &= 2q'_{nk} \log(nk/q'_{nk}) (1+o(1)) \\
        &= 2q_{n} \log(n/q_{n}) (1+o(1)).
    \end{split}
\end{equation*}
The last equality holds because we let $q'_{nk} = q_n$. 

So we have a minimax lower bound now,
$$\inf_{\{\widehat{\boldsymbol{\theta}}_{i}\}} \sup_{\{ \boldsymbol{\theta}_{i}\} \in L_0[q_n]} \sum_{i=1}^n E \lVert \widehat{\boldsymbol{\theta}}_{i} - \boldsymbol{\theta}_{i} \rVert_{\Sigma}^2 \geq 2q_{n} \log(n/q_{n}) (1+o(1)). $$
As we will see in Theorems \ref{thm3} and \ref{mEIGThm3}, when using some particular priors, the error of the Bayes estimate of $\boldsymbol{\theta}_{i}$'s will be at most $2q_{n} \log(n/q_{n}) (1+o(1))$.
This fact provides an upper bound for the minimax error, which coincides with the lower bound, and finishes this proof.

\end{proof}
\begin{proof}[Proof of Theorem \ref{thm:GL}]
We first prove \eqref{thm1}.
Observe that 
\begin{equation*}
    \lVert \widehat{\boldsymbol{\theta}}^R_i \rVert_{\Sigma}^2 = E\left(1-\kappa_i  \mid  \boldsymbol{X}_i\right)^2  \lVert\boldsymbol{X}_i\rVert_{\Sigma}^2  \leq E\left(1-\kappa_i  \mid  \boldsymbol{X}_i\right)  \lVert\boldsymbol{X}_i\rVert_{\Sigma}^2.
\end{equation*}
Denote the $l_2$-norm by $\lVert\cdot\rVert_2$.
Then, making use of Lemma \ref{TUNElemma1}, for a sequence of positive constants  $\{a_n,n\ge1\}$ to be specified later with $a_n \rightarrow \infty$ as $n \rightarrow \infty$,
\begin{equation}\label{6}
\begin{split}
    &E_0\left[ E \left(1-\kappa_i  \mid  \boldsymbol{X}_i\right)  \lVert\boldsymbol{X}_i\rVert_{\Sigma}^2   \mathbbm{1}_{\left[ \lVert\boldsymbol{X}_i\rVert_2^2  \le a_n\right]}  \right] \\
    \le& K \tau_n^a E_0\left[\exp \left(\boldsymbol{X}_i^T \boldsymbol{\Sigma}^{-1} \boldsymbol{X}_i/ 2\right)  \lVert\boldsymbol{X}_i\rVert_2^2  \mathbbm{1}_{\left[ \lVert\boldsymbol{X}_i\rVert_2^2  \le a_n\right]}  \right] \\
    =& K \tau_n^a \int_{ \lVert\boldsymbol{X}_i\rVert_2^2  \le a_n}  \lVert\boldsymbol{X}_i\rVert_2^2  \exp \left( \boldsymbol{X}_i^T \boldsymbol{\Sigma}^{-1} \boldsymbol{X}_i/2 - \boldsymbol{X}_i^T \boldsymbol{\Sigma}^{-1} \boldsymbol{X}_i / 2\right) d\boldsymbol{X}_i \\
    =& K \tau_n^a \int_{ \lVert\boldsymbol{X}_i\rVert_2^2  \le a_n}  \lVert\boldsymbol{X}_i\rVert_2^2  d\boldsymbol{X}_i \\
    \le& K \tau_n^a a_n V_k(a_n) \\
    =& K \tau_n^a a_n^{k+1},
\end{split}
\end{equation}
where $V_k(r) = \frac{\pi^{k/2}}{\Gamma(k/2+1)} r^{k}$ is the volume of a Euclidean ball of radius r in k-dimensional Euclidean space.

Moreover, under assumption (II), 
\begin{equation*}\label{7}
\begin{split}
    &E_0\left[ E \left(1-\kappa_i  \mid  \boldsymbol{X}_i\right)  \lVert\boldsymbol{X}_i\rVert_{\Sigma}^2   \mathbbm{1}_{\left[ \lVert\boldsymbol{X}_i\rVert_2^2  > a_n\right]}  \right] \\
    \le &E_0 \left( \lVert\boldsymbol{X}_i\rVert_{\Sigma}^2  \mathbbm{1}_{\left[ \lVert\boldsymbol{X}_i\rVert_{\Sigma}^2 > \lambda_{\min} \left(\Sigma^{-1}\right) a_n \right]} \right)
\end{split}
\end{equation*}

Noting that under $\boldsymbol{\theta}_{0i} = \boldsymbol{0} $, $\lVert\boldsymbol{X}_i\rVert_{\Sigma}^2 \sim \chi_k^2$, one gets
\begin{equation*}
\begin{split}
    &E_0\left(  \lVert\boldsymbol{X}_i\rVert_{\Sigma}^2 \mathbbm{1}_{\left[\lVert\boldsymbol{X}_i\rVert_{\Sigma}^2 > \lambda_{\min} \left(\Sigma^{-1}\right) a_n\right]}  \right) \\
    = & \int_{\lambda_{\min} \left(\Sigma^{-1}\right) a_n}^\infty x \exp\left(-\frac{x}{2}\right) \frac{x^{k/2-1}}{\Gamma(k/2) 2^{k/2}} dx\\
    = & 2 \int_{\left( \lambda_{\min} \left(\Sigma^{-1}\right) a_n \right)/2}^\infty \exp\left(-x\right) \frac{x^{k/2}}{\Gamma(k/2)} dx\\
    \le & 2 \exp \left(- \lambda_{\min} \left(\Sigma^{-1}\right) a_n / 4 \right) \int_{\left( \lambda_{\min} \left(\Sigma^{-1}\right) a_n \right)/2}^\infty  \exp\left(-\frac{x}{2}\right) \frac{x^{k/2}}{\Gamma(k/2)} dx\\
    \le & \exp \left(- \lambda_{\min} \left(\Sigma^{-1}\right) a_n / 4\right) \left(\frac{k}{2}\right) 2 ^{\frac{k}{2}+2}.
\end{split}
\end{equation*}

Now choosing $a_n = 4 \lambda_{\min}^{-1} \left(\Sigma^{-1}\right) (1+\epsilon_n)\log(n/q_n)$, one gets
\begin{equation}\label{8}
E_0\left[ E \left(1-\kappa_i  \mid  \boldsymbol{X}_i\right)  \lVert\boldsymbol{X}_i\rVert_{\Sigma}^2   \mathbbm{1}_{\left[ \lVert\boldsymbol{X}_i\rVert_2^2  > a_n\right]}  \right]
\le K (q_n/n)^{1+\epsilon_n}.
\end{equation}


Finally, choosing $\tau_n = \left(q_n/n\right)^{\frac{1+\epsilon_n}{a}}$, it follows from \eqref{6} and \eqref{8} that, 
\begin{equation*}
\label{9}
\begin{split}
    E_0 \left[ E\left(1-\kappa_i  \mid  \boldsymbol{X}_i\right) \lVert\boldsymbol{X}_i\rVert_{\Sigma}^2 \right] 
    & \le K \left[\left(\frac{q_n}{n}\right)^{1+\epsilon_n} \left( (1+\epsilon_n) \log \left(\frac{n}{q_n}\right)\right)^{k+1}  + \left(\frac{q_n}{n}\right)^{1+\epsilon_n} \right].
\end{split}
\end{equation*}


Summing over all $i$'s for which $\boldsymbol{\theta}_{0i} = \boldsymbol{0}$, since $\epsilon_n = 1 / \log(\log(n/q_n))$, one gets 
$$ { \sup_{\{\boldsymbol{\theta}_{0i}\} \in L_0[q_n]}} {\sum\limits_{i: \boldsymbol{\theta}_{0i} = \boldsymbol{0}} E_0 \lVert \widehat{\boldsymbol{\theta}}^R_i \rVert_{\Sigma}^2 } /
{\left( q_n \log\left(\frac{n}{q_n}\right) \right)} = o(1), \text{ as } n \rightarrow \infty.$$

Now we prove \eqref{thm2}. Use the inequality
\begin{equation*}
\begin{split}
    &E_{\boldsymbol{\theta}_{0i}} \lVert\widehat{\boldsymbol{\theta}}^R_i - \boldsymbol{\theta}_{0i}\rVert_{\Sigma}^2 \\
    =& 
    E_{\boldsymbol{\theta}_{0i}} \lVert(1-E\left(\kappa_i \mid \boldsymbol{X}_i\right))\boldsymbol{X}_i - \boldsymbol{X}_i + \boldsymbol{X}_i - \boldsymbol{\theta}_{0i}\rVert_{\Sigma}^2
    \\
    =& E_{\boldsymbol{\theta}_{0i}} \left[E\left(\kappa_i \mid \boldsymbol{X}_i\right)^2 \lVert\boldsymbol{X}_i\rVert_{\Sigma}^2 + \lVert\boldsymbol{X}_i - \boldsymbol{\theta}_{0i}\rVert_{\Sigma}^2
    -2 \langle \boldsymbol{\Sigma}^{-1/2}(\boldsymbol{X}_i - \boldsymbol{\theta}_{0i}), E\left(\kappa_i \mid \boldsymbol{X}_i\right)  \boldsymbol{\Sigma}^{-1/2}\boldsymbol{X}_i \rangle\right] \\
    \leq& E_{\boldsymbol{\theta}_{0i}} \left[E\left(\kappa_i \mid \boldsymbol{X}_i\right)^2 \lVert\boldsymbol{X}_i\rVert_{\Sigma}^2 \right] + 
    E_{\boldsymbol{\theta}_{0i}} \lVert\boldsymbol{X}_i - \boldsymbol{\theta}_{0i}\rVert_{\Sigma}^2 \\
    +& 2 
    \left( E_{\boldsymbol{\theta}_{0i}} \lVert\boldsymbol{X}_i - \boldsymbol{\theta}_{0i}\rVert_{\Sigma}^2  E_{\boldsymbol{\theta}_{0i}} \left[E\left(\kappa_i \mid \boldsymbol{X}_i\right)^2 \lVert\boldsymbol{X}_i\rVert_{\Sigma}^2 \right]  \right)^\frac{1}{2}
\end{split}
\end{equation*}

But $E_{\boldsymbol{\theta}_{0i}} \lVert\boldsymbol{X}_i - \boldsymbol{\theta}_{0i}\rVert_{\Sigma}^2 = k$, the above becomes
\begin{equation*}
\label{10}
\begin{split}
    E_{\boldsymbol{\theta}_{0i}} \lVert\widehat{\boldsymbol{\theta}}^R_i - \boldsymbol{\theta}_{0i}\rVert_{\Sigma}^2 
    \leq& E_{\boldsymbol{\theta}_{0i}} \left[E\left(\kappa_i \mid \boldsymbol{X}_i\right)^2 \lVert\boldsymbol{X}_i\rVert_{\Sigma}^2 \right] \\
    +& 2 
    k^{\frac{1}{2}} E_{\boldsymbol{\theta}_{0i}}^\frac{1}{2}  \left[E\left(\kappa_i \mid \boldsymbol{X}_i\right)^2 \lVert\boldsymbol{X}_i\rVert_{\Sigma}^2 \right] + k \\
    \leq& E_{\boldsymbol{\theta}_{0i}} \left[E\left(\kappa_i \mid \boldsymbol{X}_i\right) \lVert\boldsymbol{X}_i\rVert_{\Sigma}^2 \right] \\
    +& 2 
    k^{\frac{1}{2}} E_{\boldsymbol{\theta}_{0i}}^\frac{1}{2}  \left[E\left(\kappa_i \mid \boldsymbol{X}_i\right) \lVert\boldsymbol{X}_i\rVert_{\Sigma}^2 \right] + k
\end{split}
\end{equation*}

Since $q_n k / q_n \log(n/q_n) \rightarrow 0$, as $n \rightarrow \infty$,
it suffices to show that
\begin{equation}\label{thm2.eq1}
    \limsup\limits_{n\rightarrow\infty} { \sup_{\{\boldsymbol{\theta}_{0i}\} \in L_0[q_n]}}
    \frac{\sum\limits_{i: \boldsymbol{\theta}_{0i} \neq \boldsymbol{0}} E_{\boldsymbol{\theta}_{0i}} [E\left(\kappa_i \mid \boldsymbol{X}_i\right) \lVert\boldsymbol{X}_i\rVert_{\Sigma}^2]} 
    {2 q_n \log(n/q_n)} \leq 1.
\end{equation}
and 
\begin{equation}\label{thm2.eq2}
    \limsup\limits_{n\rightarrow\infty} { \sup_{\{\boldsymbol{\theta}_{0i}\} \in L_0[q_n]}}
    \frac{\sum\limits_{i: \boldsymbol{\theta}_{0i} \neq \boldsymbol{0}} E_{\boldsymbol{\theta}_{0i}}^{\frac{1}{2}} [E\left(\kappa_i \mid \boldsymbol{X}_i\right) \lVert\boldsymbol{X}_i\rVert_{\Sigma}^2]} 
    {q_n \log(n/q_n)} = 0.
\end{equation}

In view of Lemma \ref{TUNElemma2}, for sufficiently large $b_n > 0$, uniformly in $\boldsymbol{\theta}_{0i} \neq \boldsymbol{0}$,
\begin{equation}\label{16}
\begin{split}
    &E \left(\kappa_i \mathbbm{1}_{\left[ \kappa_i > \xi \right]} \mid  \boldsymbol{X}_i\right)
     \lVert\boldsymbol{X}_i\rVert_{\Sigma}^2 \mathbbm{1}_{\left[ \lVert\boldsymbol{X}_i\rVert_{\Sigma}^2  > b_n \right]} \\
    \leq& K \tau_n^{-a} \lVert\boldsymbol{X}_i\rVert_{\Sigma}^2
    \exp{\left[ -\frac{\xi(1-\delta)}{2} \lVert\boldsymbol{X}_i\rVert_{\Sigma}^2  \right]} \mathbbm{1}_{\left[ \lVert\boldsymbol{X}_i\rVert_{\Sigma}^2 > b_n \right]} \\ 
    \leq& K \tau_n^{-a} b_n \exp{\left[ -\frac{\xi(1-\delta)}{2} b_n \right]} \\
    \leq& K(q_n/n)^{\rho_n-\epsilon_n} \log\left(\frac{n}{q_n}\right),
\end{split}
\end{equation}
by choosing $b_n = \frac{2}{\xi(1-\delta)}(1+\rho_n) \log\left(\frac{n}{q_n}\right)$ with $\rho_n = 1/\log(\log(\log(n/q_n))$, and recalling that $\tau_n = \left(q_n/n\right)^{\frac{1+\epsilon_n}{a}}$ with $\epsilon_n = 1 / \log(\log(n/q_n))$.

Then summing over all $i$'s for which $\boldsymbol{\theta}_{0i} \neq \boldsymbol{0}$, one gets
\begin{equation}\label{21}
    \limsup\limits_{n \rightarrow \infty} { \sup_{\{\boldsymbol{\theta}_{0i}\} \in L_0[q_n]}}
    \frac{\sum\limits_{i: \boldsymbol{\theta}_{0i} \neq \boldsymbol{0}}  
    E_{\boldsymbol{\theta}_{0i}} \left[ 
    E \left(\kappa_i \mathbbm{1}_{\left[ \kappa_i > \xi \right]} \mid  \boldsymbol{X}_i\right)
     \lVert\boldsymbol{X}_i\rVert_{\Sigma}^2 \mathbbm{1}_{\left[\lVert\boldsymbol{X}_i\rVert_{\Sigma}^2 > b_n \right]}\right]}
    {q_n \log(n/q_n)} = 0.
\end{equation}

Finally,
\begin{equation}\label{17}
\begin{split}
    &E_{\boldsymbol{\theta}_{0i}} \left[ E\left(\kappa_i \mathbbm{1}_{\left[ \kappa_i > \xi \right]}  \mid  \boldsymbol{X}_i\right) \lVert\boldsymbol{X}_i\rVert_{\Sigma}^2 \mathbbm{1}_{\left[ \lVert\boldsymbol{X}_i\rVert_{\Sigma}^2 \le b_n \right]} \right] \\
    \le &
    E_{\boldsymbol{\theta}_{0i}} \left[\lVert\boldsymbol{X}_i\rVert_{\Sigma}^2 \mathbbm{1}_{\left[  \lVert\boldsymbol{X}_i\rVert_{\Sigma}^2 \le b_n \right]} \right] \\
    \le & b_n \\
    \le & \frac{2(1+\rho_n)}{\xi(1-\delta)} \log\left(\frac{n}{q_n}\right).
\end{split}
\end{equation}

Since the result above is independent of any specific choice of the parameters, by making $\xi \rightarrow 1$ and $\delta \rightarrow 0$, 
and summing over all $i$'s for which $\boldsymbol{\theta}_{0i} \neq \boldsymbol{0}$,
one gets
\begin{equation}
\begin{split}
    \limsup\limits_{n\rightarrow\infty} { \sup_{\{\boldsymbol{\theta}_{0i}\} \in L_0[q_n]}}
    \frac{\sum\limits_{i: \boldsymbol{\theta}_{0i} \neq \boldsymbol{0}} E_{\boldsymbol{\theta}_{0i}} \left[ E\left(\kappa_i \mathbbm{1}_{\left[ \kappa_i > \xi \right]}  \mid  \boldsymbol{X}_i\right)  \lVert\boldsymbol{X}_i\rVert_{\Sigma}^2 \mathbbm{1}_{\left[\lVert\boldsymbol{X}_i\rVert_{\Sigma}^2 \le b_n \right]} \right]}
    {2 q_n \log(n/q_n)} \leq 1.
\end{split}
\end{equation}
Together with Lemma \ref{TUNElemma3} and \eqref{21}, this leads to \eqref{thm2.eq1}.

Altogether Lemma \ref{TUNElemma3}, \eqref{16} and \eqref{17} also imply that
\begin{equation}
    E_{\boldsymbol{\theta}_{0i}} [E\left(\kappa_i \mid \boldsymbol{X}_i\right) \lVert\boldsymbol{X}_i\rVert_{\Sigma}^2]
    \leq K \log\left(\frac{n}{q_n}\right)(1+o(1)), \text{ as } n \rightarrow \infty.
\end{equation}

Again summing over all $i$'s for which $\boldsymbol{\theta}_{0i} \neq \boldsymbol{0}$, one gets \eqref{thm2.eq2}.
This completes the proof of \eqref{thm2}.

\eqref{thm3} follows \eqref{thm1} and \eqref{thm2}, immediately.
In particular, Theorem \ref{thm:minimax.error} provides a lower bound for the estimation error and \eqref{thm3.second} follows.
\end{proof}

\begin{proof}[Proof of Theorem \ref{GL_contraction}]
By Markov's inequality and the independence of samples,
    \begin{equation*}
        \begin{split}
            &E_{\{\boldsymbol{\theta}_{0i}\}} \Pi(\sum_{i=1}^{n} \lVert \boldsymbol{\theta}_i - \widehat{\boldsymbol{\theta}}^R_i  \rVert_{\Sigma}^2 > q_n \log(\frac{n}{q_n}) \mid \{\boldsymbol{X}_i\}) \\
            \leq &E_{\{\boldsymbol{\theta}_{0i}\}} E(\sum_{i=1}^{n} \lVert \boldsymbol{\theta}_i - \widehat{\boldsymbol{\theta}}^R_i  \rVert_{\Sigma}^2 \mid \{\boldsymbol{X}_i\}) / q_n \log(\frac{n}{q_n}) \\
            = & \sum_{i=1}^{n} E_{\boldsymbol{\theta}_{0i}} E( \lVert \boldsymbol{\theta}_i - \widehat{\boldsymbol{\theta}}^R_i  \rVert_{\Sigma}^2 \mid \boldsymbol{X}_i) / q_n \log(\frac{n}{q_n}).
        \end{split}
    \end{equation*}
Since
$$\boldsymbol{\theta}_{i} \mid \boldsymbol{X}_i, \kappa_i \sim N_k(\widehat{\boldsymbol{\theta}}^R_i, (1-\kappa_i)\boldsymbol{\Sigma}),$$
we have
    \begin{equation*}
        \begin{split}
            &E( \lVert \boldsymbol{\theta}_i - \widehat{\boldsymbol{\theta}}^R_i  \rVert_{\Sigma}^2 \mid \boldsymbol{X}_i) \\
            =& E\{ E( \lVert \boldsymbol{\theta}_i - \widehat{\boldsymbol{\theta}}^R_i  \rVert_{\Sigma}^2 \mid \boldsymbol{X}_i, \kappa_i) \mid \boldsymbol{X}_i \}\\
            = & E\{ (1-\kappa_i) E( \frac{ \lVert \boldsymbol{\theta}_i - \widehat{\boldsymbol{\theta}}^R_i  \rVert_{\Sigma}^2 }{1-\kappa_i} \mid \boldsymbol{X}_i, \kappa_i) \mid \boldsymbol{X}_i \} \\
            = & k  E (1-\kappa_i  \mid \boldsymbol{X}_i ).
        \end{split}
    \end{equation*}
    
Now we only need to find a bound for $E_{\boldsymbol{\theta}_{0i}} E (1-\kappa_i  \mid \boldsymbol{X}_i )$.
When $\boldsymbol{\theta}_{0i} \neq 0$, $$E_{\boldsymbol{\theta}_{0i}} E (1-\kappa_i  \mid \boldsymbol{X}_i ) \leq 1.$$
When $\boldsymbol{\theta}_{0i} = 0$, letting $\{a_n\}$ be a sequence of positive numbers that will be chosen later, using Lemma \ref{TUNElemma1}, 
    \begin{equation*}
        \begin{split}
            &E_{\boldsymbol{\theta}_{0i}} E (1-\kappa_i  \mid \boldsymbol{X}_i ) \\
            =& E_{\boldsymbol{0}} \{E (1-\kappa_i  \mid \boldsymbol{X}_i ) \mathbbm{1}_{\left[ \lVert\boldsymbol{X}_i\rVert_2^2  \le a_n\right]}\} + E_{\boldsymbol{0}} \{E (1-\kappa_i  \mid \boldsymbol{X}_i ) \mathbbm{1}_{\left[ \lVert\boldsymbol{X}_i\rVert_2^2  > a_n\right]}\} \\
            \leq& E_{\boldsymbol{0}} \{ K \tau_n^a \exp \left( \boldsymbol{X}_i^T \boldsymbol{\Sigma}^{-1} \boldsymbol{X}_i / 2 \right) \mathbbm{1}_{\left[ \lVert\boldsymbol{X}_i\rVert_2^2  \le a_n\right]}\} + E_{\boldsymbol{0}} \{\mathbbm{1}_{\left[ \lVert\boldsymbol{X}_i\rVert_2^2  > a_n\right]}\} \\
            =&  K \tau_n^a V_k(a_n) + K Pr (\chi^2_k > a_n) \\
            \leq &  K \tau_n^a a_n^k + K a_n^{k/2} \exp(-a_n/2),
        \end{split}
    \end{equation*}
where $V_k(r)$ is the volume of a Euclidean ball of radius r in k-dimensional Euclidean space, and the probability is bounded using the Chernoff bound for the $\chi^2$ random variables.

For $\epsilon_n = 1/\log(\log(n/q_n))$, choose $a_n = 2(1+\epsilon_n) \log (n/q_n)$ and $\tau_n = (q_n/n)^{(1+\epsilon_n)/a}$. 
Then,
    \begin{equation*}
        \begin{split}
            &\sup_{\{\boldsymbol{\theta}_{0i}\} \in L_0[q_n]} \sum_{i=1}^{n} E_{\boldsymbol{\theta}_{0i}} E( \lVert \boldsymbol{\theta}_i - \widehat{\boldsymbol{\theta}}^R_i  \rVert_{\Sigma}^2 \mid \boldsymbol{X}_i) / q_n \log(\frac{n}{q_n}) \\
            \leq &\frac{q_n + (n-q_n)(K \tau_n^a a_n^k + K a_n^{k/2} \exp(-a_n/2))}{q_n \log(\frac{n}{q_n})} \rightarrow 0,
        \end{split}
    \end{equation*}
as $n \rightarrow \infty$. 
This proves \eqref{contraction_est}.

Next, since 
$$E( \lVert \boldsymbol{\theta}_i - {\boldsymbol{\theta}}_{0i}  \rVert_{\Sigma}^2 \mid \boldsymbol{X}_i) \leq 2E( \lVert \boldsymbol{\theta}_i - \widehat{\boldsymbol{\theta}}^R_i  \rVert_{\Sigma}^2 \mid \boldsymbol{X}_i) + 2E( \lVert \widehat{\boldsymbol{\theta}}^R_i - {\boldsymbol{\theta}}_{0i}   \rVert_{\Sigma}^2 \mid \boldsymbol{X}_i),$$
by Markov's inequality, \eqref{thm3} and \eqref{contraction_est} immediately lead to \eqref{contraction_truth}.

\end{proof}

\begin{proof}[Proof of Theorem \ref{thm:mEIG}]
The proof of \eqref{mEIGThm1} is similar to the proof of \eqref{thm1}, but now using Lemma \ref{mEIGlemma1}.
We start from
\begin{equation*}
    \lVert \widehat{\boldsymbol{\theta}}^{EIG}_i \rVert_{\Sigma}^2 = E\left(1-\kappa_i  \mid  \boldsymbol{X}_i\right)^2  \lVert\boldsymbol{X}_i\rVert_{\Sigma}^2  \leq E\left(1-\kappa_i  \mid  \boldsymbol{X}_i\right)  \lVert\boldsymbol{X}_i\rVert_{\Sigma}^2.
\end{equation*}
Then, for a sequence of positive constants  $\{a_n,n\ge1\}$ to be specified later with $a_n \rightarrow \infty$ as $n \rightarrow \infty$,
\begin{equation*}
\begin{split}
    &E_0\left[ E \left(1-\kappa_i  \mid  \boldsymbol{X}_i\right)  \lVert\boldsymbol{X}_i\rVert_{\Sigma}^2   \mathbbm{1}_{\left[ \lVert\boldsymbol{X}_i\rVert_2^2  \le a_n\right]}  \right] \\
    \le& K c_n^d E_0\left[\exp \left( \boldsymbol{X}_i^T \boldsymbol{\Sigma}^{-1} \boldsymbol{X}_i /2 \right)  \lVert\boldsymbol{X}_i\rVert_{\Sigma}^2  \mathbbm{1}_{\left[ \lVert\boldsymbol{X}_i\rVert_2^2  \le a_n\right]}  \right] \\
    \le& K c_n^d  a_n^{k+1},
\end{split}
\end{equation*}
and,
\begin{equation*}
\begin{split}
    &E_0\left[ E \left(1-\kappa_i  \mid  \boldsymbol{X}_i\right)  \lVert\boldsymbol{X}_i\rVert_{\Sigma}^2   \mathbbm{1}_{\left[ \lVert\boldsymbol{X}_i\rVert_2^2  > a_n\right]}  \right] \\
    \le &E_0 \left( \lVert\boldsymbol{X}_i\rVert_{\Sigma}^2 \mathbbm{1}_{\left[ \lVert\boldsymbol{X}_i\rVert_{\Sigma}^2 > \lambda_{\min} \left(\Sigma^{-1}\right) a_n \right]} \right) \\
    \le & \exp \left(-\lambda_{\min} \left(\Sigma^{-1}\right) a_n/4\right).
\end{split}
\end{equation*}
By choosing $a_n = 4 \lambda_{\min}^{-1} \left(\Sigma^{-1}\right) (1+\epsilon_n)\log(n/q_n)$ and $c_n = \left(q_n/n\right)^{\frac{1+\epsilon_n}{d}}$, 
with $\epsilon_n = 1/\log(\log(n/q_n))$, one gets
\begin{equation*}
\begin{split}
    E_0 \left[ E\left(1-\kappa_i  \mid  \boldsymbol{X}_i\right) \lVert\boldsymbol{X}_i\rVert_{\Sigma}^2 \right] 
    & \le K \left[\left(\frac{q_n}{n}\right)^{1+\epsilon_n} \left( (1+\epsilon_n) \log \left(\frac{n}{q_n}\right)\right)^{k+1}  + \left(\frac{q_n}{n}\right)^{1+\epsilon_n} \right],
\end{split}
\end{equation*}
and hence
$${ \sup_{\{\boldsymbol{\theta}_{0i}\} \in L_0[q_n]}} {\sum\limits_{i: \boldsymbol{\theta}_{0i} = \boldsymbol{0}} E_0 \lVert \widehat{\boldsymbol{\theta}}^{EIG}_i \rVert_{\Sigma}^2} /
{\left( q_n \log\left(\frac{n}{q_n}\right) \right)} = o(1), \text{ as } n \rightarrow \infty.$$

The proof of \eqref{mEIGThm2} is similar to the proof of \eqref{thm2}, but now using Lemmas \ref{mEIGlemma2} and \ref{mEIGlemma3}.
Again, it suffices to show that
\begin{equation}\label{part1}
    \limsup\limits_{n\rightarrow\infty} { \sup_{\{\boldsymbol{\theta}_{0i}\} \in L_0[q_n]}}
    \frac{\sum\limits_{i: \boldsymbol{\theta}_{0i} \neq \boldsymbol{0}} E_{\boldsymbol{\theta}_{0i}} [E\left(\kappa_i \mid \boldsymbol{X}_i\right) \lVert\boldsymbol{X}_i\rVert_{\Sigma}^2]} 
    {2 q_n \log(n/q_n)} \leq 1.
\end{equation}
and 
\begin{equation}\label{part2}
    \limsup\limits_{n\rightarrow\infty} { \sup_{\{\boldsymbol{\theta}_{0i}\} \in L_0[q_n]}}
    \frac{\sum\limits_{i: \boldsymbol{\theta}_{0i} \neq \boldsymbol{0}} E_{\boldsymbol{\theta}_{0i}}^{\frac{1}{2}} [E\left(\kappa_i \mid \boldsymbol{X}_i\right) \lVert\boldsymbol{X}_i\rVert_{\Sigma}^2]} 
    {q_n \log(n/q_n)} = 0.
\end{equation}
By Lemma \ref{mEIGlemma2}, for sufficiently large $b_n > 0$, uniformly in $\boldsymbol{\theta}_{0i} \neq \boldsymbol{0}$,
\begin{equation*}
\begin{split}
    &E \left(\kappa_i \mathbbm{1}_{\left[ \kappa_i > \xi \right]} \mid  \boldsymbol{X}_i\right)
     \lVert\boldsymbol{X}_i\rVert_{\Sigma}^2 \mathbbm{1}_{\left[  \lVert\boldsymbol{X}_i\rVert_{\Sigma}^2 > b_n \right]} \\
    \leq& K c_n^{-d} \lVert\boldsymbol{X}_i\rVert_{\Sigma}^2
    \exp{\left[ -\frac{\xi(1-\delta)}{2} \lVert\boldsymbol{X}_i\rVert_{\Sigma}^2  \right]} \mathbbm{1}_{\left[\lVert\boldsymbol{X}_i\rVert_{\Sigma}^2 > b_n \right]} \\ 
    \leq& K c_n^{-d}  b_n \exp{\left[ -\frac{\xi(1-\delta)}{2} b_n \right]} \\
    \leq& K(q_n/n)^{\rho_n-\epsilon_n} \log\left(\frac{n}{q_n}\right),
\end{split}
\end{equation*}
by choosing $b_n = \frac{2(1+\rho_n)}{\xi(1-\delta)} \log\left(\frac{n}{q_n}\right)$, 
with $\rho = 1/ \log(\log(\log(n/q_n)))$, and recalling that $c_n = \left(q_n/n\right)^{\frac{1+\epsilon_n}{d}}$.
So,
\begin{equation}\label{part1.1}
    \limsup\limits_{n \rightarrow \infty} { \sup_{\{\boldsymbol{\theta}_{0i}\} \in L_0[q_n]}}
    \frac{\sum\limits_{i: \boldsymbol{\theta}_{0i} \neq \boldsymbol{0}}  
    E_{\boldsymbol{\theta}_{0i}} \left[ 
    E \left(\kappa_i \mathbbm{1}_{\left[ \kappa_i > \xi \right]} \mid  \boldsymbol{X}_i\right)
     \lVert\boldsymbol{X}_i\rVert_{\Sigma}^2 \mathbbm{1}_{\left[ \lVert\boldsymbol{X}_i\rVert_{\Sigma}^2 > b_n \right]}\right]}
    {q_n \log(n/q_n)} = 0.
\end{equation}
Also,
\begin{equation}\label{part1.2}
\begin{split}
    &E_{\boldsymbol{\theta}_{0i}} \left[ E\left(\kappa_i \mathbbm{1}_{\left[ \kappa_i > \xi \right]}  \mid  \boldsymbol{X}_i\right) \lVert\boldsymbol{X}_i\rVert_{\Sigma}^2 \mathbbm{1}_{\left[ \lVert\boldsymbol{X}_i\rVert_{\Sigma}^2 \le b_n \right]} \right] \\
    \le &
    E_{\boldsymbol{\theta}_{0i}} \left[\lVert\boldsymbol{X}_i\rVert_{\Sigma}^2 \mathbbm{1}_{\left[  \lVert\boldsymbol{X}_i\rVert_{\Sigma}^2 \le b_n \right]} \right] \\
    \le & b_n \\
    \le & \frac{2(1+\rho_n)}{\xi(1-\delta)} \log\left(\frac{n}{q_n}\right).
\end{split}
\end{equation}
By making $\xi \rightarrow 1$ and $\delta \rightarrow 0$,
one gets
\begin{equation}
    \limsup\limits_{n\rightarrow\infty} { \sup_{\{\boldsymbol{\theta}_{0i}\} \in L_0[q_n]}}
    \frac{\sum\limits_{i: \boldsymbol{\theta}_{0i} \neq \boldsymbol{0}} E_{\boldsymbol{\theta}_{0i}} \left[ E\left(\kappa_i \mathbbm{1}_{\left[ \kappa_i > \xi \right]}  \mid  \boldsymbol{X}_i\right)  \lVert\boldsymbol{X}_i\rVert_{\Sigma}^2 \mathbbm{1}_{\left[ \lVert\boldsymbol{X}_i\rVert_{\Sigma}^2 \le b_n \right]} \right]}
    {2 q_n \log(n/q_n)} 
    \leq 1.
\end{equation}
Together with Lemma \ref{mEIGlemma3} and \eqref{part1.1}, this leads to \eqref{part1}.

Altogether Lemma \ref{mEIGlemma3}, \eqref{part1.1} and \eqref{part1.2} also imply that
\begin{equation}
    E_{\boldsymbol{\theta}_{0i}} [E\left(\kappa_i \mid \boldsymbol{X}_i\right) \lVert\boldsymbol{X}_i\rVert_{\Sigma}^2]
    \leq K \log\left(\frac{n}{q_n}\right)(1+o(1)), \text{ as } n \rightarrow \infty.
\end{equation}
Consequently, one gets \eqref{part2}.
This completes the proof of \eqref{mEIGThm2}.

Finally, \eqref{mEIGThm3} follows the previous two results and Theorem \ref{thm:minimax.error}.
\end{proof}

\begin{proof}[Proof of Theorem \ref{mEIG_contraction}]
Similar to the Proof of Theorem \ref{GL_contraction},
we only need to find a bound for $E_{\boldsymbol{\theta}_{0i}} E (1-\kappa_i  \mid \boldsymbol{X}_i )$.
When $\boldsymbol{\theta}_{0i} \neq 0$, $$E_{\boldsymbol{\theta}_{0i}} E (1-\kappa_i  \mid \boldsymbol{X}_i ) \leq 1.$$
When $\boldsymbol{\theta}_{0i} = 0$, letting $\{a_n\}$ be a sequence of positive numbers that will be chosen later, using Lemma \ref{mEIGlemma1}, 
    \begin{equation*}
        \begin{split}
            &E_{\boldsymbol{\theta}_{0i}} E (1-\kappa_i  \mid \boldsymbol{X}_i ) \\
            =& E_{\boldsymbol{0}} \{E (1-\kappa_i  \mid \boldsymbol{X}_i ) \mathbbm{1}_{\left[ \lVert\boldsymbol{X}_i\rVert_2^2  \le a_n\right]}\} + E_{\boldsymbol{0}} \{E (1-\kappa_i  \mid \boldsymbol{X}_i ) \mathbbm{1}_{\left[ \lVert\boldsymbol{X}_i\rVert_2^2  > a_n\right]}\} \\
            \leq &  K c_n^d a_n^k + K a_n^{k/2} \exp(-a_n/2).
        \end{split}
    \end{equation*}

For some $\epsilon > 0$, choose $a_n = 2(1+\epsilon) \log (n/q_n)$ and $c_n = (q_n/n)^{(1+\epsilon)/d}$. 
Then,
    \begin{equation*}
        \begin{split}
            &\sup_{\{\boldsymbol{\theta}_{0i}\} \in L_0[q_n]} \sum_{i=1}^{n} E_{\boldsymbol{\theta}_{0i}} E( \lVert \boldsymbol{\theta}_i - \widehat{\boldsymbol{\theta}}^{EIG}_i  \rVert_{\Sigma}^2 \mid \boldsymbol{X}_i) / q_n \log(\frac{n}{q_n}) \\
            \leq &\frac{q_n + (n-q_n)(K c_n^d a_n^k + K a_n^{k/2} \exp(-a_n/2))}{q_n \log(\frac{n}{q_n})} \rightarrow 0,
        \end{split}
    \end{equation*}
as $n \rightarrow \infty$. This proves \eqref{contraction_est_mEIG}.

Next, since 
$$E( \lVert \boldsymbol{\theta}_i - {\boldsymbol{\theta}}_{0i}  \rVert_{\Sigma}^2 \mid \boldsymbol{X}_i) \leq 2E( \lVert \boldsymbol{\theta}_i - \widehat{\boldsymbol{\theta}}^{EIG}_i  \rVert_{\Sigma}^2 \mid \boldsymbol{X}_i) + 2E( \lVert \widehat{\boldsymbol{\theta}}^{EIG}_i - {\boldsymbol{\theta}}_{0i}   \rVert_{\Sigma}^2 \mid \boldsymbol{X}_i),$$
by Markov's inequality, \eqref{mEIGThm3} and \eqref{contraction_est_mEIG} immediately lead to \eqref{contraction_truth_mEIG}.

\end{proof}

\section{Proofs of Theorems in  Section \ref{sec:sets}}
\begin{proof}[Proof of Theorem \ref{CSThm}.]
We use $\widehat{\boldsymbol{\theta}}_{i}$ for $\widehat{\boldsymbol{\theta}}_{i}^R$ in this proof.

    Proof of \eqref{thm.S}. Look at the case where $\lVert \boldsymbol{\theta}_{0i} \rVert_{\Sigma}^2 \leq K_S \tau_n$.
    Note that we could write $\boldsymbol{X}_i = \boldsymbol{\theta}_{0i} + \boldsymbol{\epsilon}_i$, where $\boldsymbol{\epsilon}_i \sim \boldsymbol{N}_{k}(\boldsymbol{0}, \boldsymbol{\Sigma})$ and hence $\lVert \boldsymbol{\epsilon}_{i} \rVert_{\Sigma}^2 \sim \chi^2_k$. 
    So,
    \begin{equation*}
        \begin{split}
            \lVert \boldsymbol{\theta}_{0i} - \widehat{\boldsymbol{\theta}}_i \rVert_{\Sigma}^2 &= \lVert \boldsymbol{\theta}_{0i} - E(1-\kappa_i \mid \boldsymbol{X}_i)\boldsymbol{X}_i \rVert_{\Sigma}^2 \\
            &=\lVert E(\kappa_i \mid \boldsymbol{X}_i)\boldsymbol{\theta}_{0i} - E(1-\kappa_i \mid \boldsymbol{X}_i) \boldsymbol{\epsilon}_i \rVert_{\Sigma}^2 \\
            &\leq 2 E(\kappa_i \mid \boldsymbol{X}_i)^2 \lVert \boldsymbol{\theta}_{0i} \rVert_{\Sigma}^2 + 2 E(1-\kappa_i \mid \boldsymbol{X}_i)^2 \lVert \boldsymbol{\epsilon}_{i} \rVert_{\Sigma}^2 \\
            &\leq 2 K_s \tau_n + 2 E(1-\kappa_i \mid \boldsymbol{X}_i) \lVert \boldsymbol{\epsilon}_{i} \rVert_{\Sigma}^2.
        \end{split} 
    \end{equation*}
    By Lemma \ref{TUNElemma1}, for any fixed $\eta < a$,
    \begin{equation*}
        \begin{split}
            E(1-\kappa_i \mid \boldsymbol{X}_i) &\leq K \tau_n^a \exp \left( \frac{\boldsymbol{X}_i^T\boldsymbol{\Sigma}^{-1}\boldsymbol{X}_i}{2} \right) \\
            &\leq \tau_n^a \exp \left( \lVert \boldsymbol{\theta}_{0i} \rVert_{\Sigma}^2 \right) \exp \left( \lVert \boldsymbol{\epsilon}_{i} \rVert_{\Sigma}^2\right).
        \end{split}
    \end{equation*}
    We will show that, for small enough $\tau_n$, 
    \begin{equation}\label{LB1}
        \widehat{r}_i(\alpha, \tau_n) \geq \chi^2_{k, A\alpha} \tau_n^{1+v} (1+o(1)),
    \end{equation}
    for some $c>0$ and fixed $A > K m / (a c^a)$ with, specifically, $K = \int_0^\infty u^{-a-1} L(u) du$ here.
    Given this, for any fixed $L$, by choosing $0 < v < \rho$,
    \begin{equation*}
        \begin{split}
            P_{\boldsymbol{\theta}_{0i}} (\boldsymbol{\theta}_{0i} \in \widehat{C}^{R}_{i}) &= P_{\boldsymbol{\theta}_{0i}} (\lVert \boldsymbol{\theta}_{0i} - \widehat{\boldsymbol{\theta}}_i \rVert_{\Sigma}^2 \leq L \widehat{r}_i^{a/(1+\rho)}(\alpha, \tau_n)) \\
            &\geq P_{\boldsymbol{\theta}_{0i}} (K \tau_n + K \tau_n^a e^{\lVert \boldsymbol{\epsilon}_{i} \rVert_{\Sigma}^2} \lVert \boldsymbol{\epsilon}_{i} \rVert_{\Sigma}^2 \leq K \tau_n^{(1+v)a/(1+\rho)} (1+o(1))) \\
            &\geq P_{\boldsymbol{\theta}_{0i}} (K \tau_n + K \tau_n^a e^{\lVert \boldsymbol{\epsilon}_{i} \rVert_{\Sigma}^2} \lVert \boldsymbol{\epsilon}_{i} \rVert_{\Sigma}^2 \leq K \tau_n^{(1+v)a/(1+\rho)} (1+o(1)) \mid \lVert \boldsymbol{\epsilon}_{i} \rVert_{\Sigma}^2 \leq \chi^2_{k, \alpha}) \\
            &\times P(\lVert \boldsymbol{\epsilon}_{i} \rVert_{\Sigma}^2 \leq \chi^2_{k, \alpha})\\
            &\rightarrow 1 \times (1-\alpha) = 1-\alpha,
        \end{split}
    \end{equation*}
    as $\tau_n \rightarrow 0$, since the left hand side of the inequality in the conditional probability is of a higher order of infinitesimal.
      
    Now, it remains to prove \eqref{LB1}.
    Due to the normality of the posterior 
    \begin{equation*}
        \boldsymbol{\theta}_i \mid \boldsymbol{X}_i, \lambda_i^2 \sim \boldsymbol{N}_k ((1-\kappa_i)\boldsymbol{X}_i, (1-\kappa_i) \boldsymbol{\Sigma}),
    \end{equation*}
    applying Anderson's lemma, we have
    \begin{equation*}
        \begin{split}
            \Pi( \lVert \boldsymbol{\theta}_i - \widehat{\boldsymbol{\theta}}_{i} \rVert_{\Sigma}^2 > \widehat{r}_i(\alpha, \tau_n) \mid \boldsymbol{X}_i, \lambda_i^2) \geq \Pi( \lVert \boldsymbol{\theta}_i - (1-\kappa_i) \boldsymbol{X}_{i} \rVert_{\Sigma}^2 > \widehat{r}_i(\alpha, \tau_n) \mid \boldsymbol{X}_i, \lambda_i^2).
        \end{split}
    \end{equation*}
    Thus,
    \begin{equation*}
        \begin{split}
            \alpha &= \int_0^\infty \Pi( \lVert \boldsymbol{\theta}_i - \widehat{\boldsymbol{\theta}}_{i} \rVert_{\Sigma}^2 > \widehat{r}_i(\alpha, \tau_n) \mid \boldsymbol{X}_i, \lambda_i^2) \pi(\lambda_i^2 \mid \boldsymbol{X}_i) d\lambda_i^2 \\
            &\geq \int_0^\infty \Pi( \lVert \boldsymbol{\theta}_i - (1-\kappa_i) \boldsymbol{X}_{i} \rVert_{\Sigma}^2 > \widehat{r}_i(\alpha, \tau_n) \mid \boldsymbol{X}_i, \lambda_i^2) \pi(\lambda_i^2 \mid \boldsymbol{X}_i) d\lambda_i^2.
        \end{split}
    \end{equation*}
    Recall that $\pi(\lambda_i^2 \mid \boldsymbol{X}_i) \propto (1+\lambda_i^2 \tau_n)^{-k/2} (\lambda_i^2)^{-a-1} L(\lambda_i^2) \exp (-\frac{\boldsymbol{X}_i^T\boldsymbol{\Sigma}^{-1}\boldsymbol{X}_i}{2 (1+\lambda_i^2 \tau_n)})$.
    Let $\widetilde{\pi}(\lambda_i^2 \mid \boldsymbol{X}_i) \propto (1+\lambda_i^2 \tau_n)^{-k/2} (\lambda_i^2)^{-a-1} L(\lambda_i^2)$ be another density. Then, since $\pi(\lambda_i^2 \mid \boldsymbol{X}_i)/\widetilde{\pi}(\lambda_i^2 \mid \boldsymbol{X}_i)$ and $ \Pi( \lVert \boldsymbol{\theta}_i - (1-\kappa_i) \boldsymbol{X}_{i} \rVert_{\Sigma}^2 > \widehat{r}_i(\alpha, \tau_n) \mid \boldsymbol{X}_i, \lambda_i^2)$ are both increasing in $\lambda_i^2$
    \begin{equation*}
        \begin{split}
            &\int_0^\infty \Pi( \lVert \boldsymbol{\theta}_i - (1-\kappa_i) \boldsymbol{X}_{i} \rVert_{\Sigma}^2 > \widehat{r}_i(\alpha, \tau_n) \mid \boldsymbol{X}_i, \lambda_i^2) \pi(\lambda_i^2 \mid \boldsymbol{X}_i) d\lambda_i^2 \\
            \geq &\int_0^\infty \Pi( \lVert \boldsymbol{\theta}_i - (1-\kappa_i) \boldsymbol{X}_{i} \rVert_{\Sigma}^2 > \widehat{r}_i(\alpha, \tau_n) \mid \boldsymbol{X}_i, \lambda_i^2) \widetilde{\pi}(\lambda_i^2 \mid \boldsymbol{X}_i) d\lambda_i^2.
        \end{split}
    \end{equation*}
    On the other hand, since $\lVert \boldsymbol{\theta}_i - (1-\kappa_i) \boldsymbol{X}_{i} \rVert_{\Sigma}^2 / (1-\kappa_i) \mid \boldsymbol{X}_i, \lambda_i^2 \sim \chi^2_k$,
    and  $1-\kappa_i = (\tau_n \lambda_i^2) / (1 + \tau_n \lambda_i^2) \geq \tau_n^{1+v} (1+o(1))$ when $\lambda_i^2 \geq \tau_n^v$ for some $v > 0$. Thus, for some fixed $A(>0)$ to be determined later,
    \begin{equation*}
        \begin{split}
            &\int_0^\infty \Pi( \lVert \boldsymbol{\theta}_i - (1-\kappa_i) \boldsymbol{X}_{i} \rVert_{\Sigma}^2 > \chi^2_{k, A\alpha} \tau_n^{1+v} (1+o(1)) \mid \boldsymbol{X}_i, \lambda_i^2) \widetilde{\pi}(\lambda_i^2 \mid \boldsymbol{X}_i) d\lambda_i^2 \\         
            \geq& \int_{\tau_n^v}^\infty \Pi( \lVert \boldsymbol{\theta}_i - (1-\kappa_i) \boldsymbol{X}_{i} \rVert_{\Sigma}^2 > \chi^2_{k, A\alpha} \tau_n^{1+v} (1+o(1)) \mid \boldsymbol{X}_i, \lambda_i^2) \widetilde{\pi}(\lambda_i^2 \mid \boldsymbol{X}_i) d\lambda_i^2 \\
            \geq& \int_{\tau_n^v}^\infty \Pi( \lVert \boldsymbol{\theta}_i - (1-\kappa_i) \boldsymbol{X}_{i} \rVert_{\Sigma}^2 > \chi^2_{k, A\alpha} (1-\kappa_i) \mid \boldsymbol{X}_i, \lambda_i^2) \widetilde{\pi}(\lambda_i^2 \mid \boldsymbol{X}_i) d\lambda_i^2 \\
            =& A \alpha \widetilde{\Pi}(\lambda_i^2 \geq \tau_n^v \mid \boldsymbol{X}_i).
        \end{split}
    \end{equation*}
    Now fix $c>0$, when $n$ is sufficiently large, $\tau_n^v < c$. Since, by the dominated convergence theorem,
    \begin{equation*}
        \begin{split}
            \widetilde{\Pi}(\lambda_i^2 \geq \tau_n^v \mid \boldsymbol{X}_i) 
            = &\frac{\int_{\tau_n^v}^{\infty} (1+\lambda_i^2 \tau_n)^{-k/2} (\lambda_i^2)^{-a-1} L(\lambda_i^2) d\lambda_i^2} {\int_{0}^{\infty} (1+\lambda_i^2 \tau_n)^{-k/2} (\lambda_i^2)^{-a-1} L(\lambda_i^2) d\lambda_i^2} \\
            \geq &\frac{\int_{c}^{\infty} (1+\lambda_i^2 \tau_n)^{-k/2} (\lambda_i^2)^{-a-1} L(\lambda_i^2) d\lambda_i^2} {\int_{0}^{\infty} (1+\lambda_i^2 \tau_n)^{-k/2} (\lambda_i^2)^{-a-1} L(\lambda_i^2) d\lambda_i^2} \\
            \rightarrow &\frac{\int_{c}^{\infty} (\lambda_i^2)^{-a-1} L(\lambda_i^2) d\lambda_i^2} {\int_{0}^{\infty} (\lambda_i^2)^{-a-1} L(\lambda_i^2) d\lambda_i^2} \\
            \geq & K m \int_{c}^\infty (\lambda_i^2)^{-a-1} d \lambda_i^2 = K m / (a c^a).
        \end{split}
    \end{equation*}
    If we fix $A > a c^a/(K m)$ with $K = \int_0^\infty u^{-a-1} L(u) du$, then
    \begin{equation*}
        \begin{split}
            &\int_0^\infty \Pi( \lVert \boldsymbol{\theta}_i - (1-\kappa_i) \boldsymbol{X}_{i} \rVert_{\Sigma}^2 > \chi^2_{k, A\alpha} \tau_n^{1+v} (1+o(1)) \mid \boldsymbol{X}_i, \lambda_i^2) \widetilde{\pi}(\lambda_i^2 \mid \boldsymbol{X}_i) d\lambda_i^2 \\  
            \geq & \alpha A K m / (a c^a) \\
            > & \alpha \\
            \geq &\int_0^\infty \Pi( \lVert \boldsymbol{\theta}_i - (1-\kappa_i) \boldsymbol{X}_{i} \rVert_{\Sigma}^2 > \widehat{r}_i(\alpha, \tau_n) \mid \boldsymbol{X}_i, \lambda_i^2) \widetilde{\pi}(\lambda_i^2 \mid \boldsymbol{X}_i) d\lambda_i^2.
        \end{split}
    \end{equation*}
    This implies that $\widehat{r}_i(\alpha, \tau_n) \geq \chi^2_{k, A\alpha} \tau_n^{1+v} (1+o(1))$
    which completes the proof of the first part.

    Proof of \eqref{thm.M}
    For the case where $f_{\tau_n} \tau_n \leq \lVert \boldsymbol{\theta}_{0i} \rVert_{\Sigma}^2 \leq K_M \log \frac{1}{\tau_n}$, we start with the inequality
    \begin{equation*}
        \lVert \boldsymbol{X}_{i} \rVert_{\Sigma} \leq  \lVert \boldsymbol{\theta}_{0i} \rVert_{\Sigma} + \lVert \boldsymbol{\epsilon}_{0i} \rVert_{\Sigma}.
    \end{equation*}
    For $K_0 > K_M$,
    \begin{equation*}
        \begin{split}
            \lVert \boldsymbol{X}_{i} \rVert_{\Sigma} - \left(K_0 \log \frac{1}{\tau_n} \right)^{1/2} \leq \lVert \boldsymbol{\epsilon}_{0i} \rVert_{\Sigma} + \left(\sqrt{K_M} - \sqrt{K_0}\right) \left( \log \frac{1}{\tau_n} \right)^{1/2} \leq 0,
        \end{split}
    \end{equation*}
    if
    $$\lVert \boldsymbol{\epsilon}_{0i} \rVert_{\Sigma} \leq \left(\sqrt{K_0} - \sqrt{K_M}\right) \left( \log \frac{1}{\tau_n} \right)^{1/2},$$
    the probability of which converges to 1 as $\tau_n \rightarrow 0$.

    We now study $\widehat{r}_i(\alpha, \tau_n)$ in this case. 
    We will find an upper bound for $\widehat{r}_i(\alpha, \tau_n)$.
    As the first ingredient, for $B > 0$,
    \begin{equation*}
        \begin{split}
            &{\Pi}(\lambda_i^2 \geq B \mid \boldsymbol{X}_i) \\
            = &\frac{\int_{B}^{\infty} (1+\lambda_i^2 \tau_n)^{-k/2} (\lambda_i^2)^{-a-1} L(\lambda_i^2) \exp \left(-\frac{\lVert \boldsymbol{X}_i \rVert_{\Sigma}^2}{2 (1+\lambda_i^2 \tau_n)}\right) d\lambda_i^2} {\int_{0}^{\infty} (1+\lambda_i^2 \tau_n)^{-k/2} (\lambda_i^2)^{-a-1} L(\lambda_i^2) \exp \left(-\frac{\lVert \boldsymbol{X}_i \rVert_{\Sigma}^2}{2 (1+\lambda_i^2 \tau_n)}\right) d\lambda_i^2}.
        \end{split}
    \end{equation*}
    The denominator
    \begin{equation*}
        \begin{split}
            &\int_{0}^{\infty} (1+\lambda_i^2 \tau_n)^{-k/2} (\lambda_i^2)^{-a-1} L(\lambda_i^2) \exp \left(-\frac{\lVert \boldsymbol{X}_i \rVert_{\Sigma}^2}{2 (1+\lambda_i^2 \tau_n)}\right) d\lambda_i^2 \\
            \geq & \exp \left(-\frac{\lVert \boldsymbol{X}_i \rVert_{\Sigma}^2}{2}\right) \int_{1}^{2} (1+\lambda_i^2 \tau_n)^{-k/2} (\lambda_i^2)^{-a-1} L(\lambda_i^2)  d\lambda_i^2 \\
            \geq & \exp \left(-\frac{\lVert \boldsymbol{X}_i \rVert_{\Sigma}^2}{2}\right) (3)^{-k/2} L(1)  \int_{1}^{2}  (\lambda_i^2)^{-a-1}  d\lambda_i^2\\
            =& K \exp \left(-\frac{\lVert \boldsymbol{X}_i \rVert_{\Sigma}^2}{2}\right).
        \end{split}
    \end{equation*}
    Hence, 
    \begin{equation}
    \label{tailB}
        \begin{split}
            {\Pi}(\lambda_i^2 \geq B \mid \boldsymbol{X}_i) 
            \leq & K \exp \left(\frac{\lVert \boldsymbol{X}_i \rVert_{\Sigma}^2}{2}\right) \int_{B}^{\infty} (\lambda_i^2)^{-a-1}  d\lambda_i^2 \\
            = & K \exp \left(\frac{\lVert \boldsymbol{X}_i \rVert_{\Sigma}^2}{2}\right)B^{-a}.
        \end{split}
    \end{equation}
    Next, by the inequality
    \begin{equation*}
        \lVert \boldsymbol{\theta}_i - \widehat{\boldsymbol{\theta}}_{i} \rVert_{\Sigma}^2 \leq 2 \lVert \boldsymbol{\theta}_i - (1-\kappa_i) \boldsymbol{X}_i \rVert_{\Sigma}^2 + 2 \lVert \widehat{\boldsymbol{\theta}}_{i} - (1-\kappa_i) \boldsymbol{X}_i \rVert_{\Sigma}^2,
    \end{equation*}
    for $r > 0$,
    \begin{equation*}
        \frac{1}{2} \lVert \boldsymbol{\theta}_i - \widehat{\boldsymbol{\theta}}_{i} \rVert_{\Sigma}^2 - \lVert \widehat{\boldsymbol{\theta}}_{i} - (1-\kappa_i) \boldsymbol{X}_i \rVert_{\Sigma}^2 \geq r
    \end{equation*}
    would imply that
    \begin{equation*}
         \lVert \boldsymbol{\theta}_i - (1-\kappa_i) \boldsymbol{X}_i \rVert_{\Sigma}^2 \geq r.
    \end{equation*}
    Now,
    \begin{equation}\label{bound.reason}
        \begin{split}
            &\Pi( \frac{1}{2} \lVert \boldsymbol{\theta}_i - \widehat{\boldsymbol{\theta}}_{i} \rVert_{\Sigma}^2 - \lVert \widehat{\boldsymbol{\theta}}_{i} - (1-\kappa_i) \boldsymbol{X}_i \rVert_{\Sigma}^2 \geq r \mid \boldsymbol{X}_i, \lambda_i^2 ) \\
            \leq & \Pi( \lVert \boldsymbol{\theta}_i - (1-\kappa_i) \boldsymbol{X}_i \rVert_{\Sigma}^2 \geq r \mid \boldsymbol{X}_i, \lambda_i^2 ).
        \end{split}
    \end{equation}
    Using \eqref{bound.reason},
    \begin{equation*}
        \begin{split}
            &\int_0^\infty \Pi( \lVert \boldsymbol{\theta}_i - \widehat{\boldsymbol{\theta}}_{i} \rVert_{\Sigma}^2 \geq 2r + 2 \sup_{\lambda_i^2 \leq B} \lVert \widehat{\boldsymbol{\theta}}_{i} - (1-\kappa_i) \boldsymbol{X}_i \rVert_{\Sigma}^2 \mid \boldsymbol{X}_i, \lambda_i^2 ) \pi(\lambda_i^2 \mid \boldsymbol{X}_i) d \lambda_i^2 \\
            \leq& \int_0^B \Pi( \lVert \boldsymbol{\theta}_i - \widehat{\boldsymbol{\theta}}_{i} \rVert_{\Sigma}^2 \geq 2r + 2 \sup_{\lambda_i^2 \leq B} \lVert \widehat{\boldsymbol{\theta}}_{i} - (1-\kappa_i) \boldsymbol{X}_i \rVert_{\Sigma}^2 \mid \boldsymbol{X}_i, \lambda_i^2 ) \pi(\lambda_i^2 \mid \boldsymbol{X}_i) d \lambda_i^2 \\
            + & {\Pi}(\lambda_i^2 \geq B \mid \boldsymbol{X}_i)
            \\
            \leq& \int_0^B \Pi( \lVert \boldsymbol{\theta}_i - \widehat{\boldsymbol{\theta}}_{i} \rVert_{\Sigma}^2 \geq 2r + 2  \lVert \widehat{\boldsymbol{\theta}}_{i} - (1-\kappa_i) \boldsymbol{X}_i \rVert_{\Sigma}^2 \mid \boldsymbol{X}_i, \lambda_i^2 ) \pi(\lambda_i^2 \mid \boldsymbol{X}_i) d \lambda_i^2 \\
            + & {\Pi}(\lambda_i^2 \geq B \mid \boldsymbol{X}_i)
            \\
            \leq& \int_0^B \Pi( \lVert \boldsymbol{\theta}_i - (1-\kappa_i) \boldsymbol{X}_i \rVert_{\Sigma}^2 \geq r \mid \boldsymbol{X}_i, \lambda_i^2 ) \pi(\lambda_i^2 \mid \boldsymbol{X}_i) d \lambda_i^2 
            +  {\Pi}(\lambda_i^2 \geq B \mid \boldsymbol{X}_i).
        \end{split} 
    \end{equation*}
    Since $1-\kappa_i \leq \frac{B\tau_n}{1+B\tau_n}$ when $\lambda_i^2 \leq B$, we can bound the first term above by $\alpha/2$ via
    \begin{equation*}
        \begin{split}
            &\int_0^B \Pi( \lVert \boldsymbol{\theta}_i - (1-\kappa_i) \boldsymbol{X}_i \rVert_{\Sigma}^2 \geq \chi^2_{k,\alpha/2} \frac{B\tau_n}{1+B\tau_n} \mid \boldsymbol{X}_i, \lambda_i^2 ) \pi(\lambda_i^2 \mid \boldsymbol{X}_i) d \lambda_i^2 \\
            \leq &\int_0^B \Pi( \lVert \boldsymbol{\theta}_i - (1-\kappa_i) \boldsymbol{X}_i \rVert_{\Sigma}^2 \geq \chi^2_{k,\alpha/2} (1-\kappa_i) \mid \boldsymbol{X}_i, \lambda_i^2 ) \pi(\lambda_i^2 \mid \boldsymbol{X}_i) d \lambda_i^2 \\
            \leq& \alpha/2.
        \end{split}
    \end{equation*}
    As for the second term, by \eqref{tailB}, when $\boldsymbol{X}_i$ is fixed, for large enough $B$, 
    \begin{equation}\label{tailB2}
        {\Pi}(\lambda_i^2 \geq B \mid \boldsymbol{X}_i) \leq \alpha / 2.
    \end{equation}
    This leads to
    \begin{equation*}
        \begin{split}
            &\int_0^\infty \Pi( \lVert \boldsymbol{\theta}_i - \widehat{\boldsymbol{\theta}}_{i} \rVert_{\Sigma}^2 > \widehat{r}_i(\alpha, \tau_n) \mid \boldsymbol{X}_i, \lambda_i^2) \pi(\lambda_i^2 \mid \boldsymbol{X}_i) d\lambda_i^2\\
            = & \alpha \\
            \geq & \int_0^\infty \Pi( \lVert \boldsymbol{\theta}_i - \widehat{\boldsymbol{\theta}}_{i} \rVert_{\Sigma}^2 \geq 2r + 2 \sup_{\lambda_i^2 \leq B} \lVert \widehat{\boldsymbol{\theta}}_{i} - (1-\kappa_i) \boldsymbol{X}_i \rVert_{\Sigma}^2  \mid \boldsymbol{X}_i, \lambda_i^2 ) \pi(\lambda_i^2 \mid \boldsymbol{X}_i) d \lambda_i^2,
        \end{split} 
    \end{equation*}
    if $r = \chi^2_{k,\alpha/2} B\tau_n / (1+B\tau_n)$ and $B$ is large enough.
    Thus, when 
    $\lVert \boldsymbol{X}_{i} \rVert_{\Sigma}^2 \leq K_0 \log \frac{1}{\tau_n}$,
    for small enough $\tau_n$,
    \begin{equation*}
        \begin{split}
            \widehat{r}_i(\alpha, \tau_n) \leq& 2\chi^2_{k,\alpha/2} \frac{B\tau_n}{1+B\tau_n} + 2 \sup_{\lambda_i^2 \leq B} \lVert \widehat{\boldsymbol{\theta}}_{i} - (1-\kappa_i) \boldsymbol{X}_i \rVert_{\Sigma}^2 \\
            \leq& 2\chi^2_{k,\alpha/2} \frac{B\tau_n}{1+B\tau_n} + 4 \sup_{\lambda_i^2 \leq B} \lVert E(1-\kappa_i \mid \boldsymbol{X}_i)  {\boldsymbol{X}}_{i}  \rVert_{\Sigma}^2 +  4 \sup_{\lambda_i^2 \leq B} \lVert (1-\kappa_i) \boldsymbol{X}_i \rVert_{\Sigma}^2 \\
            \leq& 2\chi^2_{k,\alpha/2} B\tau_n + 8 B \tau_n \lVert  \boldsymbol{X}_i \rVert_{\Sigma}^2 \\
            \leq& K B \tau_n \log \frac{1}{\tau_n}.
        \end{split}
    \end{equation*}

    Again, applying Lemma \ref{TUNElemma1}, when 
    $\lVert \boldsymbol{X}_{i} \rVert_{\Sigma}^2 \leq K_0 \log \frac{1}{\tau_n}$, for fixed $\eta < a$,
    \begin{equation*}
        \begin{split}
            \lVert \widehat{\boldsymbol{\theta}}_{i}\rVert_{\Sigma}^2 \leq&
            E(1-\kappa_i \mid \boldsymbol{X}_i) \lVert \boldsymbol{X}_{i} \rVert_{\Sigma}^2 \\
            \leq& K \tau_n^{\eta} \exp \left( \lVert \boldsymbol{X}_{i} \rVert_{\Sigma}^2 / 2 \right) \lVert \boldsymbol{X}_{i} \rVert_{\Sigma}^2 \\
            \leq& K \tau_n^{\eta} \exp \left( K_0 \log \frac{1}{\tau_n} / 2 \right) K_0 \log \frac{1}{\tau_n}\\
            =& K \tau_n^{\eta - K_0 / 2} \log \frac{1}{\tau_n}.
        \end{split}
    \end{equation*}
    If we choose $B$ to be such that $B^a = K\tau_n^{-K_0/2}$, in which the factor $K$ makes \eqref{tailB2} hold, we then have, as $\tau_n \rightarrow 0$,
    \begin{equation*}
        \begin{split}
            2 \lVert \widehat{\boldsymbol{\theta}}_{i} \rVert_{\Sigma}^2 + 2 \widehat{r}_i^{a/(1+\rho)}(\alpha, \tau_n) \leq K \tau_n^{\eta - K_0 / 2} \log \frac{1}{\tau_n} + K \left( \tau_n^{1 - K_0/2a} \log \frac{1}{\tau_n} \right)^{a/(1+\rho)} = o(1),
        \end{split}
    \end{equation*}
    if we require $K_M < 2a $ and fix $\eta$ and $K_0$ to be such that $K_M < K_0 < 2\eta < 2a$. This leads to that, when 
    $\lVert \boldsymbol{X}_{i} \rVert_{\Sigma}^2 \leq K_0 \log \frac{1}{\tau_n}$, 
    for some $f_{\tau_n}$ such that $f_{\tau_n} \rightarrow \infty$ and $f_{\tau_n} \tau_n \rightarrow 0$ as $\tau_n \rightarrow 0$,
    \begin{equation*}
        \lVert {\boldsymbol{\theta}}_{0i} \rVert_{\Sigma}^2 \geq f_{\tau_n} \tau_n \geq 2 \lVert \widehat{\boldsymbol{\theta}}_{i} \rVert_{\Sigma}^2 + 2 \widehat{r}_i^{a/(1+\rho)}(\alpha, \tau_n).
    \end{equation*}
    Finally, by the inequality 
    \begin{equation*}
        \lVert {\boldsymbol{\theta}}_{0i} \rVert_{\Sigma}^2 \leq 2 \lVert {\boldsymbol{\theta}}_{0i} - \widehat{\boldsymbol{\theta}}_{i} \rVert_{\Sigma}^2 + 2 \lVert \widehat{\boldsymbol{\theta}}_{i} \rVert_{\Sigma}^2,
    \end{equation*}
    the fact that
    \begin{equation*}
        \lVert {\boldsymbol{\theta}}_{0i} \rVert_{\Sigma}^2 \geq 2 \lVert \widehat{\boldsymbol{\theta}}_{i} \rVert_{\Sigma}^2 + 2 \widehat{r}_i^{a/(1+\rho)}(\alpha, \tau_n),
    \end{equation*}
    would imply that
    \begin{equation*}
        \lVert {\boldsymbol{\theta}}_{0i} - \widehat{\boldsymbol{\theta}}_{i} \rVert_{\Sigma}^2 \geq \widehat{r}_i^{a/(1+\rho)}(\alpha, \tau_n).
    \end{equation*}
    So,
    \begin{equation*}
        \begin{split}
            &P_{\boldsymbol{\theta}_{0i}} (\boldsymbol{\theta}_{0i} \notin \widehat{C}^{R}_{i}) \\
            =& P_{\boldsymbol{\theta}_{0i}} (\lVert \boldsymbol{\theta}_{0i} - \widehat{\boldsymbol{\theta}}_i \rVert_{\Sigma}^2 > L \widehat{r}_i^{a/(1+\rho)}(\alpha, \tau_n)) \\
            \geq& P_{\boldsymbol{\theta}_{0i}} (\lVert \boldsymbol{\theta}_{0i} - \widehat{\boldsymbol{\theta}}_i \rVert_{\Sigma}^2 > L \widehat{r}_i^{a/(1+\rho)}(\alpha, \tau_n) \mid \lVert \boldsymbol{\epsilon}_{0i} \rVert_{\Sigma}^2 \leq (\sqrt{K_0} - \sqrt{K_M})^2 \log 1/\tau_n ) \\
            \times& P(\lVert \boldsymbol{\epsilon}_{0i} \rVert_{\Sigma}^2 \leq (\sqrt{K_0} - \sqrt{K_M})^2 \log 1/\tau_n)\\
            =& P(\lVert \boldsymbol{\epsilon}_{0i} \rVert_{\Sigma}^2 \leq (\sqrt{K_0} - \sqrt{K_M})^2 \log 1/\tau_n)\\
            \rightarrow& 1,
        \end{split}
    \end{equation*}
    as $\tau_n \rightarrow 0$, for any fixed $L > 0$.

    Proof of \eqref{thm.L}. Consider the case where $\lVert \boldsymbol{\theta}_{0i} \rVert_{\Sigma}^2 \geq K_L \log \frac{1}{\tau_n}$. We write
    \begin{equation*}
        \begin{split}
            \lVert \boldsymbol{\theta}_{0i} - \widehat{\boldsymbol{\theta}}_i \rVert_{\Sigma} &= \lVert \boldsymbol{\theta}_{0i} - \boldsymbol{X}_i - E(\kappa_i \mid \boldsymbol{X}_i)\boldsymbol{X}_i \rVert_{\Sigma} \\
            &\leq  \lVert \boldsymbol{\epsilon}_{i} \rVert_{\Sigma} + \lVert E(\kappa_i \mid \boldsymbol{X}_i)\boldsymbol{X}_i \rVert_{\Sigma}.
        \end{split}
    \end{equation*}
    Applying both Lemmas \ref{TUNElemma2} and \ref{TUNElemma3}, for any fixed $\xi, \delta \in (0,1)$,
    \begin{equation*}
        \begin{split}
            \lVert E(\kappa_i \mid \boldsymbol{X}_i)\boldsymbol{X}_i \rVert_{\Sigma}^2
            &\leq 2 \left(\frac{K}{\lVert \boldsymbol{X}_i \rVert_{\Sigma}^2}\right)^2 \lVert \boldsymbol{X}_i \rVert_{\Sigma}^2 \\
            &+ 2 \left( K \tau_n^{-a} \exp \left(-\frac{\xi (1-\delta)}{2} \lVert \boldsymbol{X}_i \rVert_{\Sigma}^2 \right) \right)^2 \lVert \boldsymbol{X}_i \rVert_{\Sigma}^2,
        \end{split}
    \end{equation*}
    Further, 
    $$\lVert \boldsymbol{X}_i \rVert_{\Sigma} = \lVert \boldsymbol{\theta}_{0i} + \boldsymbol{\epsilon}_i \rVert_{\Sigma} \geq \lvert \lVert \boldsymbol{\theta}_{0i} \rVert_{\Sigma} - \lVert \boldsymbol{\epsilon}_i \rVert_{\Sigma} \rvert,$$
    in which $\lVert \boldsymbol{\theta}_{0i} \rVert_{\Sigma}^2 \geq  K_L \log \frac{1}{\tau_n}$, by assumption, and $\lVert \boldsymbol{\epsilon}_i \rVert_{\Sigma}^2 \leq \chi^2_{k, \alpha}$ with probability $1-\alpha$.
    Since $\tau_n \rightarrow 0$, when $\lVert \boldsymbol{\epsilon}_i \rVert_{\Sigma}^2 \leq \chi^2_{k, \alpha}$, we will have
    $$\lVert \boldsymbol{X}_i \rVert_{\Sigma} \geq \left( K_L \log \frac{1}{\tau_n} \right)^{1/2} - \left( \chi^2_{k, \alpha}\right)^{1/2}$$
    and, consequently,
    $$\lVert \boldsymbol{X}_i \rVert_{\Sigma}^2 \geq K_L \log \frac{1}{\tau_n} (1+o(1)).$$
    By choosing $K_L$, $\xi$, and $\delta$ to be such that $K_L  (1+o(1)) > \frac{2a}{\xi (1-\delta)}$, e.g., $K_L = 3a$ and $\xi(1-\delta) = 3/4$, when $\lVert \boldsymbol{\epsilon}_i \rVert_{\Sigma}^2 \leq \chi^2_{k, \alpha}$, we will have, as $\tau_n \rightarrow 0$,
    \begin{equation*}
            \lVert E(\kappa_i \mid \boldsymbol{X}_i)\boldsymbol{X}_i \rVert_{\Sigma}^2 \leq o(1),
    \end{equation*}
    hence,
    \begin{equation}\label{diff.bound}
            \begin{split}
            \lVert \boldsymbol{\theta}_{0i} - \widehat{\boldsymbol{\theta}}_i \rVert_{\Sigma}^2 &\leq \left( \sqrt{\chi^2_{k,\alpha}} + o(1) \right)^2 = \chi^2_{k,\alpha} (1 + o(1)).
        \end{split}
    \end{equation}
    Then, we find a lower bound for $\widehat{r}_i(\alpha, \tau_n)$ in this case.
    Making use of the posterior normality and Anderson's lemma again,
    \begin{equation*}
        \begin{split}
            \alpha &= \int_0^\infty \Pi( \lVert \boldsymbol{\theta}_i - \widehat{\boldsymbol{\theta}}_{i} \rVert_{\Sigma}^2 > \widehat{r}_i(\alpha, \tau_n) \mid \boldsymbol{X}_i, \lambda_i^2) \pi(\lambda_i^2 \mid \boldsymbol{X}_i) d\lambda_i^2 \\
            &\geq \int_0^\infty \Pi( \lVert \boldsymbol{\theta}_i - (1-\kappa_i) \boldsymbol{X}_{i} \rVert_{\Sigma}^2 > \widehat{r}_i(\alpha, \tau_n) \mid \boldsymbol{X}_i, \lambda_i^2) \pi(\lambda_i^2 \mid \boldsymbol{X}_i) d\lambda_i^2.
        \end{split}
    \end{equation*}
    On the other hand, since $1-\kappa_i \geq g_{\tau_n} / (1+g_{\tau_n}) = 1 + o(1),$ if $\lambda_i^2 \geq g_{\tau_n}/\tau_n$ with $g_{\tau_n} = (\log 1/\tau_n)^{1/3}$, for some fixed $A$,
    \begin{equation*}
        \begin{split}
            &\int_0^\infty \Pi( \lVert \boldsymbol{\theta}_i - (1-\kappa_i) \boldsymbol{X}_{i} \rVert_{\Sigma}^2 > \chi^2_{k, A\alpha} g_{\tau_n} / (1+g_{\tau_n}) \mid \boldsymbol{X}_i, \lambda_i^2) {\pi}(\lambda_i^2 \mid \boldsymbol{X}_i) d\lambda_i^2 \\         
            \geq& \int_{g_{\tau_n}/\tau_n}^\infty \Pi( \lVert \boldsymbol{\theta}_i - (1-\kappa_i) \boldsymbol{X}_{i} \rVert_{\Sigma}^2 > \chi^2_{k, A\alpha} g_{\tau_n} / (1+g_{\tau_n}) \mid \boldsymbol{X}_i, \lambda_i^2) {\pi}(\lambda_i^2 \mid \boldsymbol{X}_i) d\lambda_i^2 \\
            \geq& \int_{g_{\tau_n}/\tau_n}^\infty \Pi( \lVert \boldsymbol{\theta}_i - (1-\kappa_i) \boldsymbol{X}_{i} \rVert_{\Sigma}^2 > \chi^2_{k, A\alpha} (1-\kappa_i) \mid \boldsymbol{X}_i, \lambda_i^2) {\pi}(\lambda_i^2 \mid \boldsymbol{X}_i) d\lambda_i^2 \\
            =& A \alpha {\Pi}(\lambda_i^2 \geq g_{\tau_n}/\tau_n \mid \boldsymbol{X}_i).
        \end{split}
    \end{equation*}
    We now study the posterior probability in a situation where $\lVert \boldsymbol{X}_i \rVert_{\Sigma}^2 \geq K_L \log \frac{1}{\tau_n} (1+o(1))$:
    \begin{equation*}
        \begin{split}
            &{\Pi}(\lambda_i^2 < g_{\tau_n}/\tau_n \mid \boldsymbol{X}_i) \\
            = &\frac{\int^{g_{\tau_n}/\tau_n}_{0} (1+\lambda_i^2 \tau_n)^{-k/2} (\lambda_i^2)^{-a-1} L(\lambda_i^2) \exp \left(-\frac{\lVert \boldsymbol{X}_i \rVert_{\Sigma}^2}{2 (1+\lambda_i^2 \tau_n)}\right) d\lambda_i^2} {\int_{0}^{\infty} (1+\lambda_i^2 \tau_n)^{-k/2} (\lambda_i^2)^{-a-1} L(\lambda_i^2) \exp \left(-\frac{\lVert \boldsymbol{X}_i \rVert_{\Sigma}^2}{2 (1+\lambda_i^2 \tau_n)}\right) d\lambda_i^2} \\
            := & N/D.
        \end{split}
    \end{equation*}
    Next,
    \begin{equation*}
        \begin{split}
            D \geq &\int_{(g_{\tau_n}+1)/\tau_n}^{2(g_{\tau_n}+1)/\tau_n} (1+\lambda_i^2 \tau_n)^{-k/2} (\lambda_i^2)^{-a-1} L(\lambda_i^2) \exp \left(-\frac{\lVert \boldsymbol{X}_i \rVert_{\Sigma}^2}{2 (1+\lambda_i^2 \tau_n)}\right) d\lambda_i^2 \\
            \geq & m (2 g_{\tau_n}+3)^{-k/2} \exp \left(-\frac{\lVert \boldsymbol{X}_i \rVert_{\Sigma}^2}{2 (2+g_{\tau_n})}\right) \int_{(g_{\tau_n}+1)/\tau_n}^{2(g_{\tau_n}+1)/\tau_n} (\lambda_i^2)^{-a-1} d\lambda_i^2 \\
            =& K (2 g_{\tau_n}+3)^{-k/2} \exp \left(-\frac{\lVert \boldsymbol{X}_i \rVert_{\Sigma}^2}{2 (2+g_{\tau_n})}\right) \left(\frac{\tau_n}{g_{\tau_n}+1}\right)^a.
        \end{split}
    \end{equation*}
    Next, fixing a constant $c > 0$, when $g_{\tau_n}$ is large enough,
    \begin{equation*}
        \begin{split}
            N = &\int^{g_{\tau_n}/\tau_n}_{0} (1+\lambda_i^2 \tau_n)^{-k/2} (\lambda_i^2)^{-a-1} L(\lambda_i^2) \exp \left(-\frac{\lVert \boldsymbol{X}_i \rVert_{\Sigma}^2}{2 (1+\lambda_i^2 \tau_n)}\right) d\lambda_i^2\\
            = & \left(\int^{c/\tau_n}_{0} + \int^{g_{\tau_n}/\tau_n}_{c/\tau_n} \right)(1+\lambda_i^2 \tau_n)^{-k/2} (\lambda_i^2)^{-a-1} L(\lambda_i^2) \exp \left(-\frac{\lVert \boldsymbol{X}_i \rVert_{\Sigma}^2}{2 (1+\lambda_i^2 \tau_n)}\right) d\lambda_i^2 \\
            =: & N_1 + N_2.
        \end{split}
    \end{equation*}
    The first term
    \begin{equation*}
        \begin{split}
            N_1 = &\int^{c/\tau_n}_{0} (1+\lambda_i^2 \tau_n)^{-k/2} (\lambda_i^2)^{-a-1} L(\lambda_i^2) \exp \left(-\frac{\lVert \boldsymbol{X}_i \rVert_{\Sigma}^2}{2 (1+\lambda_i^2 \tau_n)}\right) d\lambda_i^2 \\
            \leq &\int^{c/\tau_n}_{0} (\lambda_i^2)^{-a-1} L(\lambda_i^2) \exp \left(-\frac{\lVert \boldsymbol{X}_i \rVert_{\Sigma}^2}{2 (1+c)}\right) d\lambda_i^2 \\
            \leq &  \exp \left(-\frac{\lVert \boldsymbol{X}_i \rVert_{\Sigma}^2}{2 (1+c)}\right) \int^{\infty}_{0} (\lambda_i^2)^{-a-1} L(\lambda_i^2) d\lambda_i^2 \\
            = & K \exp \left(-\frac{\lVert \boldsymbol{X}_i \rVert_{\Sigma}^2}{2 (1+c)}\right) .
        \end{split}
    \end{equation*}
    When $\lVert \boldsymbol{X}_i \rVert_{\Sigma}^2 \geq K_L \log \frac{1}{\tau_n} (1+o(1))$,
    \begin{equation}\label{N1}
        \begin{split}
            N_1/D \leq & K (2 g_{\tau_n}+3)^{k/2} \left(\frac{g_{\tau_n}+1}{\tau_n}\right)^a \exp \left(\frac{\lVert \boldsymbol{X}_i \rVert_{\Sigma}^2}{2 (2+g_{\tau_n})}-\frac{\lVert \boldsymbol{X}_i \rVert_{\Sigma}^2}{2 (1+c)}\right) \\
            \leq& K (2 g_{\tau_n}+3)^{k/2} \left(\frac{g_{\tau_n}+1}{\tau_n}\right)^a \exp \left(-\frac{K_L}{2} \log \frac{1}{\tau_n} \frac{1}{(1+c)} (1+o(1)) \right) \\
            =& K (2 g_{\tau_n}+3)^{k/2} \left(g_{\tau_n}+1\right)^a \tau_n^{\frac{K_L}{2(1+c)}(1+o(1))-a} = o(1),
        \end{split}
    \end{equation}
    as $\tau_n \rightarrow 0$, if $c$ is chosen properly, e.g., $c = 1/3$. Our choice of $K_L$ makes sure that the above is $o(1)$ as $\tau_n \rightarrow 0$, given that $g_{\tau_n} = (\log 1/\tau_n)^{1/3}$.
    Moreover,
    \begin{equation*}
        \begin{split}
            N_2 = &\int_{c/\tau_n}^{g_{\tau_n}/\tau_n} (1+\lambda_i^2 \tau_n)^{-k/2} (\lambda_i^2)^{-a-1} L(\lambda_i^2) \exp \left(-\frac{\lVert \boldsymbol{X}_i \rVert_{\Sigma}^2}{2 (1+\lambda_i^2 \tau_n)}\right) d\lambda_i^2 \\
            \leq &M (1+c)^{-k/2} \exp \left(-\frac{\lVert \boldsymbol{X}_i \rVert_{\Sigma}^2}{2 (1+g_{\tau_n})}\right) \int_{c/\tau_n}^{g_{\tau_n}/\tau_n} (\lambda_i^2)^{-a-1} d\lambda_i^2 \\
            \leq &M (1+c)^{-k/2} \exp \left(-\frac{\lVert \boldsymbol{X}_i \rVert_{\Sigma}^2}{2 (1+g_{\tau_n})}\right) (\tau_n / c)^{a+1} \frac{g_{\tau_n} - c}{\tau_n}.
        \end{split}
    \end{equation*}
    It follows that 
    \begin{equation*}
        \begin{split}
            N_2/D \leq & K (2 g_{\tau_n}+3)^{k/2} \left(g_{\tau_n}+1\right)^a (g_{\tau_n} - c) \exp \left(-\frac{\lVert \boldsymbol{X}_i \rVert_{\Sigma}^2}{2 (1+g_{\tau_n})(2+g_{\tau_n})}\right).
        \end{split}
    \end{equation*}
    Since we have chosen $g_{\tau_n} = (\log 1/\tau_n)^{1/3}$, when $\lVert \boldsymbol{X}_i \rVert_{\Sigma}^2 \geq K_L \log \frac{1}{\tau_n} (1+o(1))$, we will have, as $\tau_n \rightarrow 0$,
    \begin{equation}\label{N2}
        N_2/D = o(1).
    \end{equation} 
    Combining \eqref{N1} and \eqref{N2}, when $\lVert \boldsymbol{X}_i \rVert_{\Sigma}^2 \geq K_L \log \frac{1}{\tau_n} (1+o(1))$, as $\tau_n \rightarrow 0$,
    \begin{equation*}
        \begin{split}
            {\Pi}(\lambda_i^2 < g_{\tau_n}/\tau_n \mid \boldsymbol{X}_i) = o(1),
        \end{split}
    \end{equation*}
    and hence
    \begin{equation*}
        \begin{split}
            &\int_0^\infty \Pi( \lVert \boldsymbol{\theta}_i - (1-\kappa_i) \boldsymbol{X}_{i} \rVert_{\Sigma}^2 > \chi^2_{k, A\alpha} g_{\tau_n} / (1+g_{\tau_n}) \mid \boldsymbol{X}_i, \lambda_i^2) {\pi}(\lambda_i^2 \mid \boldsymbol{X}_i) d\lambda_i^2 \\
            \geq& A \alpha {\Pi}(\lambda_i^2 \geq g_{\tau_n}/\tau_n \mid \boldsymbol{X}_i) \\
            =& A \alpha (1-o(1)).
        \end{split}
    \end{equation*}
    Here, let $A = \beta/\alpha$ for some fixed $\beta > \alpha$. For small enough $\tau_n$, we will have
    \begin{equation*}
        \begin{split}
            &\int_0^\infty \Pi( \lVert \boldsymbol{\theta}_i - (1-\kappa_i) \boldsymbol{X}_{i} \rVert_{\Sigma}^2 > \chi^2_{k, \beta} g_{\tau_n} / (1+g_{\tau_n}) \mid \boldsymbol{X}_i, \lambda_i^2) {\pi}(\lambda_i^2 \mid \boldsymbol{X}_i) d\lambda_i^2 \\
            \geq& \beta(1-o(1)) \\
            >& \alpha \\
            \geq & \int_0^\infty \Pi( \lVert \boldsymbol{\theta}_i - (1-\kappa_i) \boldsymbol{X}_{i} \rVert_{\Sigma}^2 > \widehat{r}_i(\alpha, \tau_n) \mid \boldsymbol{X}_i, \lambda_i^2) \pi(\lambda_i^2 \mid \boldsymbol{X}_i) d\lambda_i^2.
        \end{split}
    \end{equation*}
    It follows that 
    \begin{equation}\label{LB2}
        \widehat{r}_i(\alpha, \tau_n) \geq \chi^2_{k, \beta} g_{\tau_n} / (1+g_{\tau_n}) = \chi^2_{k, \beta} (1+o(1)),
    \end{equation} when $\lVert \boldsymbol{X}_i \rVert_{\Sigma}^2 \geq K_L \log \frac{1}{\tau_n} (1+o(1))$.
    
    As a reminder, when $\tau_n$ is small enough, $\lVert \boldsymbol{\epsilon}_i \rVert_{\Sigma}^2 \leq \chi^2_{k, \alpha}$ implies that $\lVert \boldsymbol{X}_i \rVert_{\Sigma}^2 \geq K_L \log \frac{1}{\tau_n} (1+o(1))$.
    
    So, by \eqref{diff.bound} along with \eqref{LB2}, in the case where $\lVert \boldsymbol{\theta}_{0i} \rVert_{\Sigma}^2 \geq K_L \log \frac{1}{\tau_n}$, for any fixed $\beta > \alpha$, as $\tau_n \rightarrow 0$,
    \begin{equation*}
        \begin{split}
            P_{\boldsymbol{\theta}_{0i}} (\boldsymbol{\theta}_{0i} \in \widehat{C}^{R}_{i}) &= P_{\boldsymbol{\theta}_{0i}} (\lVert \boldsymbol{\theta}_{0i} - \widehat{\boldsymbol{\theta}}_i \rVert_{\Sigma}^2 \leq L \widehat{r}_i^{a/(1+\rho)}(\alpha, \tau_n)) \\
            &\geq P_{\boldsymbol{\theta}_{0i}} (\lVert \boldsymbol{\theta}_{0i} - \widehat{\boldsymbol{\theta}}_i \rVert_{\Sigma}^2 \leq L \widehat{r}_i^{a/(1+\rho)}(\alpha, \tau_n) \mid \lVert \boldsymbol{\epsilon}_{i} \rVert_{\Sigma}^2 \leq \chi^2_{k, \alpha}) \\
            &\times P(\lVert \boldsymbol{\epsilon}_{i} \rVert_{\Sigma}^2 \leq \chi^2_{k, \alpha})\\
            &\geq P_{\boldsymbol{\theta}_{0i}} (\chi^2_{k,\alpha}(1 + o(1)) \leq L (\chi^2_{k, \beta})^{a/(1+\rho)} (1+o(1)) \mid \lVert \boldsymbol{\epsilon}_{i} \rVert_{\Sigma}^2 \leq \chi^2_{k, \alpha}) \\
            &\times P(\lVert \boldsymbol{\epsilon}_{i} \rVert_{\Sigma}^2 \leq \chi^2_{k, \alpha})\\
            &\rightarrow 1 \times (1-\alpha) = 1-\alpha,
        \end{split}
    \end{equation*}
    if we choose $L > \chi^2_{k,\alpha}(\chi^2_{k,\beta})^{-a/(1+\rho)}$, e.g., $L = 2\chi^2_{k,\alpha}(\chi^2_{k,\beta})^{-a/(1+\rho)}$.
\end{proof}

\begin{proof}[Proof of Theorem \ref{CSThm.EIG}.]
We use $\widehat{\boldsymbol{\theta}}_{i}$ for $\widehat{\boldsymbol{\theta}}_{i}^{EIG}$ in this proof.

    Proof of \eqref{thm.S.EIG}. 
    We first find a lower bound for $\widehat{r}_i(\alpha, c_n)$.
    Recall that $\pi(\kappa_i \mid \boldsymbol{X}_i) \propto \kappa_i^{d+k/2-1} (1 - \kappa_i + c_n \kappa_i)^{-d-1}  \exp (-\kappa_i \boldsymbol{X}_i^T\boldsymbol{\Sigma}^{-1}\boldsymbol{X}_i/2)$.
    Similar to the proof of \eqref{thm.S}, let $\widetilde{\pi}(\kappa_i \mid \boldsymbol{X}_i) \propto  \kappa_i^{d+k/2-1} (1 - \kappa_i + c_n \kappa_i)^{-d-1}$ be another density.
    Also, under the setup of this theorem, 
    \begin{equation*}
        \boldsymbol{\theta}_i \mid \boldsymbol{X}_i, \kappa_i \sim \boldsymbol{N}_k ((1-\kappa_i)\boldsymbol{X}_i, (1-\kappa_i) \boldsymbol{\Sigma}),
    \end{equation*}
    Now we can proceed similarly as the proof of \eqref{LB1}.
    For some $A(>0)$ and some $v (> 1)$ to be determined later,
    \begin{equation*}
        \begin{split}
            &\int_0^1 \Pi( \lVert \boldsymbol{\theta}_i - (1-\kappa_i) \boldsymbol{X}_{i} \rVert_{\Sigma}^2 > \chi^2_{k, A\alpha} c_n^v \mid \boldsymbol{X}_i, \kappa_i) \widetilde{\pi}(\kappa_i \mid \boldsymbol{X}_i) d\kappa_i \\         
            \geq& \int_0^{ 1-c_n^v} \Pi( \lVert \boldsymbol{\theta}_i - (1-\kappa_i) \boldsymbol{X}_{i} \rVert_{\Sigma}^2 > \chi^2_{k, A\alpha} c_n^v \mid \boldsymbol{X}_i, \kappa_i) \widetilde{\pi}(\kappa_i \mid \boldsymbol{X}_i) d\kappa_i \\
            \geq& \int_0^{ 1-c_n^v} \Pi( \lVert \boldsymbol{\theta}_i - (1-\kappa_i) \boldsymbol{X}_{i} \rVert_{\Sigma}^2 > \chi^2_{k, A\alpha} (1-\kappa_i) \mid \boldsymbol{X}_i, \kappa_i) \widetilde{\pi}(\kappa_i \mid \boldsymbol{X}_i) d\kappa_i \\
            =& A \alpha \widetilde{\Pi}(\kappa_i \leq  1-c_n^v \mid \boldsymbol{X}_i).
        \end{split}
    \end{equation*}
    Since
    \begin{equation*}
        \begin{split}
            \widetilde{\Pi}(\kappa_i > 1-c_n^v \mid \boldsymbol{X}_i) &= \frac{\int_{ 1-c_n^v}^1 \kappa_i^{d+k/2-1} (1 - \kappa_i + c_n \kappa_i)^{-d-1} d\kappa_i}{\int_0^1 \kappa_i^{d+k/2-1} (1 - \kappa_i + c_n \kappa_i)^{-d-1} d\kappa_i} \\
            &\leq \frac{\int_{1-c_n^v}^1 c_n^{-d-1} d\kappa_i}{\int_{1-c_n}^1 (1-c_n)^{d+k/2-1} (c_n + c_n (1-c_n))^{-d-1} d\kappa_i} \\
            &= c_n^{v-1} \frac{(2-c_n)^{d+1}}{(1-c_n)^{d+k/2-1}},
        \end{split}
    \end{equation*}
    suppose $\alpha < 1/2$, when $c_n$ is small enough, if we fix $A = 2$, for example, 
    \begin{equation*}
        A \alpha \widetilde{\Pi}(\kappa_i \leq 1-c_n^v \mid \boldsymbol{X}_i) > \alpha.
    \end{equation*}
    On the other hand, by Anderson's lemma,
    \begin{equation*}
        \begin{split}
            \alpha &= \int_0^1 \Pi( \lVert \boldsymbol{\theta}_i - \widehat{\boldsymbol{\theta}}_{i} \rVert_{\Sigma}^2 > \widehat{r}_i(\alpha, c_n) \mid \boldsymbol{X}_i, \kappa_i) \pi(\kappa_i \mid \boldsymbol{X}_i) d\kappa_i \\
            &\geq \int_0^1 \Pi( \lVert \boldsymbol{\theta}_i - (1-\kappa_i) \boldsymbol{X}_{i} \rVert_{\Sigma}^2 > \widehat{r}_i(\alpha, c_n) \mid \boldsymbol{X}_i, \kappa_i) \pi(\kappa_i \mid \boldsymbol{X}_i) d\kappa_i.
        \end{split}
    \end{equation*}
    So,
    \begin{equation*}
        \begin{split}
            &\int_0^1 \Pi( \lVert \boldsymbol{\theta}_i - (1-\kappa_i) \boldsymbol{X}_{i} \rVert_{\Sigma}^2 > \chi^2_{k, A\alpha} c_n^v \mid \boldsymbol{X}_i, \kappa_i) \widetilde{\pi}(\kappa_i \mid \boldsymbol{X}_i) d\kappa_i\\ 
            >  \alpha 
            \geq &\int_0^1 \Pi( \lVert \boldsymbol{\theta}_i - (1-\kappa_i) \boldsymbol{X}_{i} \rVert_{\Sigma}^2 > \widehat{r}_i(\alpha, c_n) \mid \boldsymbol{X}_i, \kappa_i) \widetilde{\pi}(\kappa_i \mid \boldsymbol{X}_i) d\kappa_i.
        \end{split}
    \end{equation*}
    This implies that
    \begin{equation}\label{LB1.EIG}
        \widehat{r}_i(\alpha, c_n) \geq \chi^2_{k, A\alpha} c_n^v.
    \end{equation}

    For the case where $\lVert \boldsymbol{\theta}_{0i} \rVert_{\Sigma}^2 \leq K_S' c_n$.
    Again,
    \begin{equation*}
        \begin{split}
            \lVert \boldsymbol{\theta}_{0i} - \widehat{\boldsymbol{\theta}}_i \rVert_{\Sigma}^2 &= \lVert \boldsymbol{\theta}_{0i} - E(1-\kappa_i \mid \boldsymbol{X}_i)\boldsymbol{X}_i \rVert_{\Sigma}^2 \\
            &=\lVert E(\kappa_i \mid \boldsymbol{X}_i)\boldsymbol{\theta}_{0i} - E(1-\kappa_i \mid \boldsymbol{X}_i) \boldsymbol{\epsilon}_i \rVert_{\Sigma}^2 \\
            &\leq 2 E(\kappa_i \mid \boldsymbol{X}_i)^2 \lVert \boldsymbol{\theta}_{0i} \rVert_{\Sigma}^2 + 2 E(1-\kappa_i \mid \boldsymbol{X}_i)^2 \lVert \boldsymbol{\epsilon}_{i} \rVert_{\Sigma}^2 \\
            &\leq 2 K_s \tau_n + 2 E(1-\kappa_i \mid \boldsymbol{X}_i) \lVert \boldsymbol{\epsilon}_{i} \rVert_{\Sigma}^2.
        \end{split} 
    \end{equation*}
    Using Lemma \ref{mEIGlemma1} and \eqref{LB1.EIG}, when we choose $v$ to be such that $v{d/(1+\rho)} < d < 1$, e.g., $v = 1+\rho/2$,
    \begin{equation*}
        \begin{split}
            P_{\boldsymbol{\theta}_{0i}} (\boldsymbol{\theta}_{0i} \in \widehat{C}^{EIG}_{i}) &= P_{\boldsymbol{\theta}_{0i}} (\lVert \boldsymbol{\theta}_{0i} - \widehat{\boldsymbol{\theta}}_i \rVert_{\Sigma}^2 \leq L \widehat{r}_i^{d/(1+\rho)} (\alpha, c_n)) \\
            &\geq P_{\boldsymbol{\theta}_{0i}} (K c_n + K c_n^d e^{\lVert \boldsymbol{\epsilon}_{i} \rVert_{\Sigma}^2} \lVert \boldsymbol{\epsilon}_{i} \rVert_{\Sigma}^2 \leq (\chi^2_{k, A\alpha} c_n^v)^{d/(1+\rho)} ) \\
            &\geq P_{\boldsymbol{\theta}_{0i}} (K c_n + K c_n^d e^{\lVert \boldsymbol{\epsilon}_{i} \rVert_{\Sigma}^2} \lVert \boldsymbol{\epsilon}_{i} \rVert_{\Sigma}^2 \leq (\chi^2_{k, A\alpha} c_n^v)^{d/(1+\rho)} \mid \lVert \boldsymbol{\epsilon}_{i} \rVert_{\Sigma}^2 \leq \chi^2_{k, \alpha}) \\
            &\times P(\lVert \boldsymbol{\epsilon}_{i} \rVert_{\Sigma}^2 \leq \chi^2_{k, \alpha})\\
            &\rightarrow 1 \times (1-\alpha) = 1-\alpha,
        \end{split}
    \end{equation*}
    as $c_n \rightarrow 0$, since the left hand side of the inequality in the conditional probability is of a higher order of infinitesimal.

    Proof of \eqref{thm.M.EIG}
    For the case where $ f'_{c_n} c_n \leq \lVert \boldsymbol{\theta}_{0i} \rVert_{\Sigma}^2 \leq K'_M \log \frac{1}{c_n}$, similar to the proof of \eqref{thm.M},
    for $K_0' > K_M'$,
    \begin{equation*} 
        \begin{split}
            \lVert \boldsymbol{X}_{i} \rVert_{\Sigma} - \left(K_0' \log \frac{1}{c_n} \right)^{1/2} \leq \lVert \boldsymbol{\epsilon}_{0i} \rVert_{\Sigma} + \left(\sqrt{K_M'} - \sqrt{K_0'}\right) \left( \log \frac{1}{c_n} \right)^{1/2} \leq 0,
        \end{split}
    \end{equation*}
    if
    $$\lVert \boldsymbol{\epsilon}_{0i} \rVert_{\Sigma} \leq \left(\sqrt{K_0'} - \sqrt{K_M'}\right) \left( \log \frac{1}{c_n} \right)^{1/2},$$
    the probability of which converges to 1 as $c_n \rightarrow 0$.
    
    Using \eqref{bound.reason} again,
    \begin{equation*}
        \begin{split}
            &\int_0^1 \Pi( \lVert \boldsymbol{\theta}_i - \widehat{\boldsymbol{\theta}}_{i} \rVert_{\Sigma}^2 \geq 2r + 2 \sup_{\kappa_i \geq {B_n}} \lVert \widehat{\boldsymbol{\theta}}_{i} - (1-\kappa_i) \boldsymbol{X}_i \rVert_{\Sigma}^2 \mid \boldsymbol{X}_i, \kappa_i ) \pi(\kappa_i \mid \boldsymbol{X}_i) d \kappa_i \\
            \leq& \int_{B_n}^1 \Pi( \lVert \boldsymbol{\theta}_i - \widehat{\boldsymbol{\theta}}_{i} \rVert_{\Sigma}^2 \geq 2r + 2 \sup_{\kappa_i \geq {B_n}} \lVert \widehat{\boldsymbol{\theta}}_{i} - (1-\kappa_i) \boldsymbol{X}_i \rVert_{\Sigma}^2 \mid \boldsymbol{X}_i, \kappa_i ) \pi(\kappa_i \mid \boldsymbol{X}_i) d \kappa_i \\
            + & {\Pi}(\kappa_i < {B_n} \mid \boldsymbol{X}_i)
            \\
            \leq& \int_{B_n}^1  \Pi( \lVert \boldsymbol{\theta}_i - \widehat{\boldsymbol{\theta}}_{i} \rVert_{\Sigma}^2 \geq 2r + 2  \lVert \widehat{\boldsymbol{\theta}}_{i} - (1-\kappa_i) \boldsymbol{X}_i \rVert_{\Sigma}^2 \mid \boldsymbol{X}_i, \kappa_i ) \pi(\kappa_i \mid \boldsymbol{X}_i) d \kappa_i \\
            + & {\Pi}(\kappa_i < {B_n} \mid \boldsymbol{X}_i)
            \\
            \leq& \int_{B_n}^1 \Pi( \lVert \boldsymbol{\theta}_i - (1-\kappa_i) \boldsymbol{X}_i \rVert_{\Sigma}^2 \geq r \mid \boldsymbol{X}_i, \kappa_i ) \pi(\kappa_i \mid \boldsymbol{X}_i) d \kappa_i 
            +  {\Pi}(\kappa_i < {B_n} \mid \boldsymbol{X}_i).
        \end{split} 
    \end{equation*}
    For the first term, let $r = \chi^2_{k,\alpha/2} (1-{B_n})$,
    \begin{equation*}
        \begin{split}
            &\int_{B_n}^1 \Pi( \lVert \boldsymbol{\theta}_i - (1-\kappa_i) \boldsymbol{X}_i \rVert_{\Sigma}^2 \geq \chi^2_{k,\alpha/2} (1-{B_n}) \mid \boldsymbol{X}_i, \kappa_i ) \pi(\kappa_i \mid \boldsymbol{X}_i) d \kappa_i \\
            \leq &\int_{B_n}^1 \Pi( \lVert \boldsymbol{\theta}_i - (1-\kappa_i) \boldsymbol{X}_i \rVert_{\Sigma}^2 \geq \chi^2_{k,\alpha/2} (1-\kappa_i) \mid \boldsymbol{X}_i, \kappa_i ) \pi(\kappa_i \mid \boldsymbol{X}_i) d \kappa_i \\
            \leq& \alpha/2.
        \end{split}
    \end{equation*}
    For the other term, when $\boldsymbol{X}_i$ is fixed, if ${B_n} = 1 - c_n^{d/(2d+2)}$, as $c_n \rightarrow 0$,
    \begin{equation*}
        \begin{split}
            &{\Pi}(\kappa_i < {B_n} \mid \boldsymbol{X}_i) \\
            =& \frac{\int_{0}^{{1-c_n^{d/(2d+2)}}} \kappa_i^{d+k/2-1} (1 - \kappa_i + c_n \kappa_i)^{-d-1}  \exp (-\kappa_i \boldsymbol{X}_i^T\boldsymbol{\Sigma}^{-1}\boldsymbol{X}_i/2) d \kappa_i  } {\int_{0}^{1}\kappa_i^{d+k/2-1} (1 - \kappa_i + c_n \kappa_i)^{-d-1}  \exp (-\kappa_i \boldsymbol{X}_i^T\boldsymbol{\Sigma}^{-1}\boldsymbol{X}_i/2) d \kappa_i  } \\
            \leq& \frac{\int_{0}^{{1-c_n^{d/(2d+2)}}} (c_n^{d/(2d+2)} + c_n (1-c_n^{d/(2d+2)}))^{-d-1} d \kappa_i  } {\exp(-\boldsymbol{X}_i^T\boldsymbol{\Sigma}^{-1}\boldsymbol{X}_i/2) \int_{1-c_n}^{1} (1-c_n)^{d+k/2-1} (c_n + c_n (1-c_n))^{-d-1}  d \kappa_i } \\
            =& \frac{c_n^{-d/2} (1+c_n^{(d+2)/(2d+2)}-c_n)^{-d-1}} {\exp(-\boldsymbol{X}_i^T\boldsymbol{\Sigma}^{-1}\boldsymbol{X}_i/2) c_n^{-d} (1-c_n)^{d+k/2-1} (2-c_n)^{-d-1}  }  \\ 
            \leq& K \exp(\boldsymbol{X}_i^T\boldsymbol{\Sigma}^{-1}\boldsymbol{X}_i/2) {c_n}^{d/2} \leq \alpha/2.
        \end{split}
    \end{equation*}
    So,
    \begin{equation*}
        \begin{split}
            &\int_0^1 \Pi( \lVert \boldsymbol{\theta}_i - \widehat{\boldsymbol{\theta}}_{i} \rVert_{\Sigma}^2 > \widehat{r}_i(\alpha, c_n) \mid \boldsymbol{X}_i, \kappa_i) \pi(\kappa_i \mid \boldsymbol{X}_i) d\kappa_i\\
            = & \alpha \\
            \geq& \int_0^1 \Pi( \lVert \boldsymbol{\theta}_i - \widehat{\boldsymbol{\theta}}_{i} \rVert_{\Sigma}^2 \geq 2r + 2 \sup_{\kappa_i \geq {B_n}} \lVert \widehat{\boldsymbol{\theta}}_{i} - (1-\kappa_i) \boldsymbol{X}_i \rVert_{\Sigma}^2 \mid \boldsymbol{X}_i, \kappa_i ) \pi(\kappa_i \mid \boldsymbol{X}_i) d \kappa_i,
        \end{split} 
    \end{equation*}
    if $r = \chi^2_{k,\alpha/2} (1-{B_n})$ and ${B_n} = 1 - c_n^{d/(2d+2)}$.
    Thus, when 
    $\lVert \boldsymbol{X}_{i} \rVert_{\Sigma}^2 \leq K_0' \log \frac{1}{c_n}$,
    for small enough $c_n$,
    \begin{equation*}
        \begin{split}
            \widehat{r}_i(\alpha, c_n) \leq& 2\chi^2_{k,\alpha/2} (1-{B_n}) + 2 \sup_{\kappa_i \geq {B_n}} \lVert \widehat{\boldsymbol{\theta}}_{i} - (1-\kappa_i) \boldsymbol{X}_i \rVert_{\Sigma}^2 \\
            \leq& 2\chi^2_{k,\alpha/2} c_n^{d/(2d+2)} + 4 \sup_{\kappa_i \geq {B_n}} \lVert E(1-\kappa_i \mid \boldsymbol{X}_i)  {\boldsymbol{X}}_{i}  \rVert_{\Sigma}^2 +  4 \sup_{\kappa_i \geq {B_n}} \lVert (1-\kappa_i) \boldsymbol{X}_i \rVert_{\Sigma}^2 \\
            \leq& 2\chi^2_{k,\alpha/2} c_n^{d/(2d+2)} + 8 c_n^{d/(2d+2)} \lVert \boldsymbol{X}_i \rVert_{\Sigma}^2 \\
            \leq& K c_n^{d/(2d+2)} \log \frac{1}{c_n}.
        \end{split}
    \end{equation*}

    Finally, by the inequality 
    \begin{equation*}
        \lVert {\boldsymbol{\theta}}_{0i} \rVert_{\Sigma}^2 \leq 2 \lVert {\boldsymbol{\theta}}_{0i} - \widehat{\boldsymbol{\theta}}_{i} \rVert_{\Sigma}^2 + 2 \lVert \widehat{\boldsymbol{\theta}}_{i} \rVert_{\Sigma}^2,
    \end{equation*}
    the fact that
    \begin{equation*}
        \lVert {\boldsymbol{\theta}}_{0i} \rVert_{\Sigma}^2 \geq 2 \lVert \widehat{\boldsymbol{\theta}}_{i} \rVert_{\Sigma}^2 + 2 \widehat{r}_i^{d/(1+\rho)}(\alpha, \tau_n),
    \end{equation*}
    would imply that
    \begin{equation*}
        \lVert {\boldsymbol{\theta}}_{0i} - \widehat{\boldsymbol{\theta}}_{i} \rVert_{\Sigma}^2 \geq \widehat{r}_i^{d/(1+\rho)}(\alpha, \tau_n).
    \end{equation*}
    Applying Lemma \ref{mEIGlemma1}, when 
    $\lVert \boldsymbol{X}_{i} \rVert_{\Sigma}^2 \leq K_0' \log \frac{1}{c_n}$,
    \begin{equation*}
        \begin{split}
            \lVert \widehat{\boldsymbol{\theta}}_{i}\rVert_{\Sigma}^2 \leq&
            E(1-\kappa_i \mid \boldsymbol{X}_i) \lVert \boldsymbol{X}_{i} \rVert_{\Sigma}^2 \\
            \leq& K c_n^d \exp \left( \lVert \boldsymbol{X}_{i} \rVert_{\Sigma}^2 / 2 \right) \lVert \boldsymbol{X}_{i} \rVert_{\Sigma}^2 \\
            \leq& K c_n^d \exp \left( K_0' \log \frac{1}{c_n} / 2 \right) K_0' \log \frac{1}{c_n}\\
            =& K c_n^{d - K_0' / 2} \log \frac{1}{c_n}.
        \end{split}
    \end{equation*}
    We then have, as $c_n \rightarrow 0$,
    \begin{equation*}
        \begin{split}
            2 \lVert \widehat{\boldsymbol{\theta}}_{i} \rVert_{\Sigma}^2 + 2 \widehat{r}_i^{d/(1+\rho)}(\alpha, c_n) \leq K c_n^{d - K_0' / 2} \log \frac{1}{c_n} + K \left( c_n^{d/(2d+2)} \log \frac{1}{c_n}\right)^{d/(1+\rho)} = o(1),
        \end{split}
    \end{equation*}
    if we require $K_M < 2d $ and fix $K_0'$ to be such that $K_M < K_0' < 2d$. This implies that when 
    $\lVert \boldsymbol{X}_{i} \rVert_{\Sigma}^2 \leq K_0' \log \frac{1}{c_n}$, 
    for some $f_{\tau_n}$ such that $f_{c_n}' \rightarrow \infty$ and $f_{c_n}' c_n \rightarrow 0$ as $c_n \rightarrow 0$,
    \begin{equation*}
        \lVert {\boldsymbol{\theta}}_{0i} \rVert_{\Sigma}^2 \geq f_{c_n}' c_n \geq 2 \lVert \widehat{\boldsymbol{\theta}}_{i} \rVert_{\Sigma}^2 + 2 \widehat{r}_i^{d/(1+\rho)}(\alpha, c_n).
    \end{equation*}
    Finally,    
    \begin{equation*}
        \begin{split}
            &P_{\boldsymbol{\theta}_{0i}} (\boldsymbol{\theta}_{0i} \notin \widehat{C}^{EIG}_{i}) \\
            =& P_{\boldsymbol{\theta}_{0i}} (\lVert \boldsymbol{\theta}_{0i} - \widehat{\boldsymbol{\theta}}_i \rVert_{\Sigma}^2 > L \widehat{r}_i^{d/(1+\rho)}(\alpha, \tau_n)) \\
            \geq& P_{\boldsymbol{\theta}_{0i}} (\lVert \boldsymbol{\theta}_{0i} - \widehat{\boldsymbol{\theta}}_i \rVert_{\Sigma}^2 > L \widehat{r}_i^{d/(1+\rho)}(\alpha, \tau_n) \mid \lVert \boldsymbol{\epsilon}_{0i} \rVert_{\Sigma}^2 \leq (\sqrt{K_0} - \sqrt{K_M})^2 \log 1/\tau_n ) \\
            \times& P(\lVert \boldsymbol{\epsilon}_{0i} \rVert_{\Sigma}^2 \leq (\sqrt{K_0} - \sqrt{K_M})^2 \log 1/\tau_n)\\
            =& P(\lVert \boldsymbol{\epsilon}_{0i} \rVert_{\Sigma}^2 \leq (\sqrt{K_0} - \sqrt{K_M})^2 \log 1/\tau_n)\\
            \rightarrow& 1,
        \end{split}
    \end{equation*}
    as $c_n \rightarrow 0$, for any fixed $L > 0$.

    Proof of \eqref{thm.L.EIG}. Similar to the proof of \eqref{thm.L}, using Lemmas \ref{mEIGlemma2} and \ref{mEIGlemma3}, when $\lVert \boldsymbol{\theta}_{0i} \rVert_{\Sigma}^2 \geq K_L' \log \frac{1}{c_n}$ and $\lVert \boldsymbol{\epsilon}_i \rVert_{\Sigma}^2 \leq \chi^2_{k, \alpha}$, we have, as $c_n \rightarrow 0$,
    $$\lVert \boldsymbol{X}_i \rVert_{\Sigma}^2 \geq K_L' \log \frac{1}{c_n} (1+o(1)),$$
    and, for $K_L'  (1+o(1)) > \frac{2d}{\xi (1-\delta)}$,
    \begin{equation}\label{diff.bound.EIG}
            \begin{split}
            \lVert \boldsymbol{\theta}_{0i} - \widehat{\boldsymbol{\theta}}_i \rVert_{\Sigma}^2 &\leq  \chi^2_{k,\alpha} (1 + o(1)).
        \end{split}
    \end{equation}

    Then, we find a lower bound for $\widehat{r}_i(\alpha, c_n)$ in this case.
    Making use of the posterior normality and Anderson's lemma again,
    \begin{equation*}
        \begin{split}
            \alpha &= \int_0^1 \Pi( \lVert \boldsymbol{\theta}_i - \widehat{\boldsymbol{\theta}}_{i} \rVert_{\Sigma}^2 > \widehat{r}_i(\alpha, c_n) \mid \boldsymbol{X}_i, \kappa_i) \pi(\kappa_i \mid \boldsymbol{X}_i) d\kappa_i \\
            &\geq \int_0^1 \Pi( \lVert \boldsymbol{\theta}_i - (1-\kappa_i) \boldsymbol{X}_{i} \rVert_{\Sigma}^2 > \widehat{r}_i(\alpha, c_n) \mid \boldsymbol{X}_i, \kappa_i) \pi(\kappa_i \mid \boldsymbol{X}_i) d\kappa_i.
        \end{split}
    \end{equation*}
    On the other hand, since $1-\kappa_i \geq 1 - c_n$ if $\kappa_i \leq c_n$, for some fixed $A$,
    \begin{equation*}
        \begin{split}
            &\int_0^1 \Pi( \lVert \boldsymbol{\theta}_i - (1-\kappa_i) \boldsymbol{X}_{i} \rVert_{\Sigma}^2 > \chi^2_{k, A\alpha} (1 - c_n) \mid \boldsymbol{X}_i, \kappa_i) {\pi}(\kappa_i \mid \boldsymbol{X}_i) d\kappa_i \\         
            \geq& \int_{0}^{c_n} \Pi( \lVert \boldsymbol{\theta}_i - (1-\kappa_i) \boldsymbol{X}_{i} \rVert_{\Sigma}^2 > \chi^2_{k, A\alpha} (1 - c_n) \mid \boldsymbol{X}_i, \kappa_i) {\pi}(\kappa_i \mid \boldsymbol{X}_i) d\kappa_i \\
            \geq& \int_{0}^{c_n} \Pi( \lVert \boldsymbol{\theta}_i - (1-\kappa_i) \boldsymbol{X}_{i} \rVert_{\Sigma}^2 > \chi^2_{k, A\alpha} (1-\kappa_i) \mid \boldsymbol{X}_i, \kappa_i) {\pi}(\kappa_i \mid \boldsymbol{X}_i) d\kappa_i \\
            =& A \alpha {\Pi}( \kappa_i \leq c_n \mid \boldsymbol{X}_i).
        \end{split}
    \end{equation*}
    When $\lVert \boldsymbol{X}_i \rVert_{\Sigma}^2 \geq K_L' \log \frac{1}{c_n} (1+o(1)),$ as $n \rightarrow \infty$,
    \begin{equation*}
        \begin{split}
            &{\Pi}(\kappa_i \geq c_n \mid \boldsymbol{X}_i) \\
            \leq& \frac{\int_{c_n}^{1} \kappa_i^{d+k/2-1} (1 - \kappa_i + c_n \kappa_i)^{-d-1}  \exp (-\kappa_i \boldsymbol{X}_i^T\boldsymbol{\Sigma}^{-1}\boldsymbol{X}_i/2) d\kappa_i }{\int_{1-c_n^{1/d}}^{1} \kappa_i^{d+k/2-1} (1 - \kappa_i + c_n \kappa_i)^{-d-1}  \exp (-\kappa_i \boldsymbol{X}_i^T\boldsymbol{\Sigma}^{-1}\boldsymbol{X}_i/2) d\kappa_i} \\
            \leq& \frac{\int_{c_n}^{1} c_n^{-d-1}  \exp (-c_n \boldsymbol{X}_i^T\boldsymbol{\Sigma}^{-1}\boldsymbol{X}_i/2) d\kappa_i }{\int_{1-c_n^{1/d}}^1 (1-c_n^{1/d})^{d+k/2-1} (c_n^{1/d} + c_n (1-c_n^{1/d}))^{-d-1}  \exp(-\boldsymbol{X}_i^T\boldsymbol{\Sigma}^{-1}\boldsymbol{X}_i/2) d\kappa_i} \\ 
            =& c_n^{-d} \exp(- (1-c_n) \boldsymbol{X}_i^T\boldsymbol{\Sigma}^{-1}\boldsymbol{X}_i/2) (1+o(1)) \\
            \leq& c_n^{-d} \exp \left(- (1-c_n) K_L' \log \frac{1}{c_n} (1+o(1)) /2 \right) (1+o(1)) \\
            =& c_n^{(1-c_n) K_L' (1+o(1)) /2 - d} = o(1).
        \end{split}
    \end{equation*}
    since we have already chosen $K_L'  (1+o(1)) > \frac{2d}{\xi (1-\delta)} > 2d$.

    Let $A = \beta/\alpha$ for some fixed $\beta > \alpha$. For small enough $c_n$, we will have
    \begin{equation*}
        \begin{split}
            &\int_0^1 \Pi( \lVert \boldsymbol{\theta}_i - (1-\kappa_i) \boldsymbol{X}_{i} \rVert_{\Sigma}^2 > \chi^2_{k, A\alpha} (1 - c_n) \mid \boldsymbol{X}_i, \kappa_i) {\pi}(\kappa_i \mid \boldsymbol{X}_i) d\kappa_i \\
            \geq& \beta(1-o(1)) \\
            >& \alpha \\
            \geq& \int_0^1 \Pi( \lVert \boldsymbol{\theta}_i - (1-\kappa_i) \boldsymbol{X}_{i} \rVert_{\Sigma}^2 > \widehat{r}_i(\alpha, c_n) \mid \boldsymbol{X}_i, \kappa_i) \pi(\kappa_i \mid \boldsymbol{X}_i) d\kappa_i.
        \end{split}
    \end{equation*}
    It follows that, as $c_n \rightarrow 0$,
    \begin{equation}\label{LB2.EIG}
        \widehat{r}_i(\alpha, c_n) \geq \chi^2_{k, \beta}(1 - c_n) = \chi^2_{k, \beta} (1+o(1)),
    \end{equation} when $\lVert \boldsymbol{X}_i \rVert_{\Sigma}^2 \geq K_L' \log \frac{1}{c_n} (1+o(1))$.
    
    So, by \eqref{diff.bound.EIG} along with \eqref{LB2.EIG}, in the case where $\lVert \boldsymbol{\theta}_{0i} \rVert_{\Sigma}^2 \geq K_L' \log \frac{1}{c_n}$, for any fixed $\beta > \alpha$, as $c_n \rightarrow 0$,
    \begin{equation*}
        \begin{split}
            P_{\boldsymbol{\theta}_{0i}} (\boldsymbol{\theta}_{0i} \in \widehat{C}^{EIG}_{i}) &= P_{\boldsymbol{\theta}_{0i}} (\lVert \boldsymbol{\theta}_{0i} - \widehat{\boldsymbol{\theta}}_i \rVert_{\Sigma}^2 \leq L \widehat{r}_i^{d/(1+\rho)}(\alpha, c_n)) \\
            &\geq P_{\boldsymbol{\theta}_{0i}} (\lVert \boldsymbol{\theta}_{0i} - \widehat{\boldsymbol{\theta}}_i \rVert_{\Sigma}^2 \leq L \widehat{r}_i^{d/(1+\rho)}(\alpha, c_n) \mid \lVert \boldsymbol{\epsilon}_{i} \rVert_{\Sigma}^2 \leq \chi^2_{k, \alpha}) \\
            &\times P(\lVert \boldsymbol{\epsilon}_{i} \rVert_{\Sigma}^2 \leq \chi^2_{k, \alpha})\\
            &\geq P_{\boldsymbol{\theta}_{0i}} (\chi^2_{k,\alpha}(1 + o(1)) \leq L (\chi^2_{k, \beta})^{d/(1+\rho)} (1+o(1)) \mid \lVert \boldsymbol{\epsilon}_{i} \rVert_{\Sigma}^2 \leq \chi^2_{k, \alpha}) \\
            &\times P(\lVert \boldsymbol{\epsilon}_{i} \rVert_{\Sigma}^2 \leq \chi^2_{k, \alpha})\\
            &\rightarrow 1 \times (1-\alpha) = 1-\alpha,
        \end{split}
    \end{equation*}
    if we choose $L > \chi^2_{k,\alpha}(\chi^2_{k,\beta})^{-{d/(1+\rho)}}$, e.g., $L = 2\chi^2_{k,\alpha}(\chi^2_{k,\beta})^{-{d/(1+\rho)}}$.
\end{proof}
\end{appendices}

\section*{Statements and Declarations}
\subsection*{Acknowledgement}
We thank the anonymous reviewers for their careful reading of our manuscript and their many
insightful comments and suggestions.
\subsection*{Funding}
The authors did not receive support from any organization for the submitted work.
\subsection*{Competing interests}
The authors have no competing interests to declare that are relevant to the content of this article.

 \bibliography{bibliography}


\end{document}